\documentclass[3p,sort&compress]{elsarticle}
\usepackage[software,hardware,siam]{mymacros} 
\usepackage{enumitem}
\usepackage{caption}
\usepackage{subcaption} 
\usepackage{multirow}
\usepackage{slashbox}
\usepackage{booktabs}
\usepackage[compatible]{algpseudocode} 
\usepackage[ruled]{algorithm}

\usepackage[linktocpage=true,colorlinks=true]{hyperref}

\usepackage{ulem}

\usepackage{pgfplots}

\newtheorem{theorem}{Theorem}
\newtheorem{definition}[theorem]{Definition}
\newtheorem{lemma}[theorem]{Lemma}

\newtheorem{proposition}[theorem]{Proposition}

\newtheorem{remark}[theorem]{Remark}

\usepackage[linktocpage=true,colorlinks=true]{hyperref}

\hypersetup{pdftitle={Efficient Solution of  Large-Scale Algebraic Riccati Equations 
        Associated with Index-2 DAEs via the Inexact Low-Rank Newton-ADI Method}}
\numberwithin{equation}{section}

\journal{Applied Numerical Mathematics}

\begin{document}

\begin{frontmatter}

\title{Efficient Solution of  Large-Scale Algebraic Riccati Equations 
        Associated with Index-2 DAEs via the  Inexact Low-Rank Newton-ADI Method}
        
\author[mpi,ovgu]{Peter Benner}
\ead{benner@mpi-magdeburg.mpg.de}

\author[rice]{Matthias Heinkenschloss\corref{cor1}\fnref{heinken}}
\ead{heinken@rice.edu}

\author[mpi]{Jens Saak}
\ead{saak@mpi-magdeburg.mpg.de}

\author[HKW]{Heiko~K. Weichelt\fnref{weichelt}}
\ead{heiko.weichelt@mathworks.co.uk}

\cortext[cor1]{Corresponding author}
\fntext[heinken]{The research of this author was supported in part by
   NSF grant DMS-1522798 and by the DARPA EQUiPS Program, Award UTA15-001068.}
\fntext[weichelt]{This paper is based on the PhD Thesis of this author,
  which was completed while he was with the Max Planck Institute for Dynamics of Complex Technical
  Systems Magdeburg, Germany.}
  
\address[mpi]{Research Group Computational Methods in Systems and Control
  Theory (CSC),\\ Max Planck Institute for Dynamics of Complex Technical
  Systems Magdeburg,\\ Sandtorstr. 1, 39106 Magdeburg, Germany}
\address[rice]{Department of Computational and Applied Mathematics
  (CAAM),\\ Rice University, MS-134, 6100 Main Street, Houston, TX
  77005-1892, USA}
\address[HKW]{The Mathworks Ltd., Matrix House, Cambridge Business Park,\\
CB4 0HH Cambridge, United Kingdom}  
\address[ovgu]{Institut f\"ur Analysis und Numerik, Fakult\"at f\"ur Mathematik,\\
             Otto-von-Guericke Universit\"at Magdeburg,\\
             Universit\"atsplatz 2, 39106 Magdeburg, Germany.}


\begin{abstract}
  This paper extends the algorithm of 
  {\em Benner, Heinkenschloss, Saak, and Weichelt:
  An  inexact low-rank {N}ewton-{ADI} method for large-scale algebraic {R}iccati
  equations, Applied Numerical Mathematics Vol.~108 (2016), pp.~125--142, doi:10.1016/j.apnum.2016.05.006}
  to Riccati equations associated with Hessenberg index-2 Differential Algebratic Equation (DAE) systems.
  Such DAE systems arise, e.g., from semi-discretized, linearized (around steady state) Navier-Stokes equations.
  The solution of the associated Riccati equation is important, e.g., to compute feedback laws that
  stabilize the  Navier-Stokes equations. 
  Challenges in the numerical  solution of the Riccati equation 
  arise from the large-scale of the underlying systems and the algebraic constraint in the DAE system.
  These challenges are met by a careful extension of the inexact low-rank {N}ewton-{ADI} method 
  to the case of DAE systems. A main ingredient in the extension to the DAE case is the projection
  onto the manifold described by the algebraic constraints. In the algorithm, the equations are never explicitly
  projected, but the projection is only applied as needed.  
  Numerical experience indicates that the algorithmic choices for the control of inexactness and line-search
  can help avoid subproblems with matrices that are only marginally stable.
  The performance of the algorithm is illustrated on a large-scale Riccati equation associated with the
  stabilization of Navier-Stokes flow around a cylinder.
\end{abstract}


\end{frontmatter}


\section{Introduction}
This paper introduces and analyzes an efficient algorithm for the solution of the
generalized continuous algebraic Riccati equation (GCARE) associated with 
the solution of linear quadratic regulator (LQR) problems governed by Hessenberg 
index-2 Differential Algebraic Equations (DAEs). This problem arises, e.g., in the
computation of feedback laws that stabilize Navier-Stokes flows. The numerical
solution of the Riccati equation is challenging because the underlying systems are
large-scale and because of the presence of algebraic constraints in the DAE system.
To overcome these
challenges we extend our inexact low-rank {N}ewton-{ADI} method developed
in~\cite{BenHSetal16} for problems governed by ordinary differential equations (ODEs)
to this DAE case. The main idea is to use the structure of the Hessenberg index-2 DAE
and apply the discrete version of the Leray projector 
(see Heinkenschloss et al.\ \cite{morHeiSS08} and B{\"a}nsch, et al.\ \cite{BaeBSetal15}) 
to transform the LQR problem governed by the DAE into a classical LQR problem governed 
by an ODE.\@ In principle, the standard LQR and Riccati theory as well as the inexact low-rank {N}ewton-{ADI} 
method developed in our previous paper~\cite{BenHSetal16} can be applied to this ODE problem. 
This, however
leads to a solution approach that is not practical because the projected systems
are large-scale and, because of the projection, dense. To arrive at an efficient
algorithm, the computations must be presented in terms of the original
large-scale sparse system and the structure of the governing DAE system must be
exploited. This is done in this paper.
In addition, numerical experience with our new algorithm indicates that our control of inexactness and the line-search 
leads to a start-up phase that reaches the quadratic convergence region of the Newton iteration faster and
tends to avoid marginally stable subproblems during intermediate iterations.

The LQR problem and associated Riccati equation considered in this paper have also
been solved by B{\"a}nsch et al.\ \cite{BaeBSetal15}. 
However, the focus of~\cite{BaeBSetal15} was the
computation of feedback laws for Navier-Stokes flows, and a basic version of
an inexact low-rank {N}ewton-{ADI} method was applied. Our paper focusses on
the solution of the Riccati equation and incorporates many recent improvements.
As a result, the algorithm in this paper delivers an approximately 90-times 
speed-up over the algorithm used in~\cite{BaeBSetal15}.
Benner and Stykel~\cite{BenS14} study the solution of projected Riccati equations, which are
associated with DAEs. They use so-called spectral projectors, which project onto the
right and left deflating subspaces. While these projectors can be applied to general DAEs 
defined by a regular pencil, in the general case
``the projectors [\dots] are required in explicit form [and the] computation of these projectors is, in general, 
 very expensive''~\cite[p. 590]{BenS14}. The projector used in our paper
 is specially designed for the index-2 DAE system arising for fluid flow problems and our
 Kleinman-Newton-ADI method contains many improvements not yet available in~\cite{BenS14}.
 In  principle it is possible to use rational Krylov subspace projection methods 
(see Simoncini et al.\ \cite{Sim16,SimSM14}) to solve the Riccati equations, but extensions
of this approach to the DAE case and numerical comparisons of the latest versions of both
approaches are not yet available.

As pointed out above, a main ingredient for the efficiency of our approach is the 
exploitation of the special structure of the  Hessenberg index-2 DAE, in what is
called implicit index-reduction. Specifically, we can use structured projectors, rather
than generic and expensive spectral projectors.  Implicit index-reduction can also
be applied to other structured DAE systems, see e.g.~\cite{morFreRM08,morGugSW13,morSaaV18,morBenSU16a}. 
We demonstrate our approach on a large-scale Riccati equation associated with the stabilization of 
Navier-Stokes flow, but the extension of the techniques described in this paper to other
saddle point structured DAEs is straight forward.

This paper is organized as follows. 
The next section, Section~\ref{sec:riccati}, introduces the LQR problem, uses projection onto the constraint
manifold to derive a projected Riccati equation, and reviews existence results
for both the projected Riccati equation and the LQR problem.
Section~\ref{sec:newt_kleinm} reviews the main components of our
algorithm in~\cite{BenHSetal16} applied to the projected GCARE and
Section~\ref{sec:newt_kleinm_DAE2} carefully exploits the  special structure of
the projected GCARE 
for an efficient numerical realization of the inexact low-rank {N}ewton-{ADI} method.
Finally, Section~\ref{sec:numerics} illustrates the performance of our algorithm 
on a large-scale Riccati equation associated with the stabilization of Navier-Stokes 
flow around a cylinder --- a problem also solved by B{\"a}nsch et al.\ \cite{BaeBSetal15}. 
As mentioned earlier, the algorithmic improvements in this paper lead to approximately 90-times 
speed-up over the algorithm used in~\cite{BaeBSetal15}.

{\bf Notation.}  Throughout the paper we consider the Hilbert space of
matrices in $\mathbb{R}^{n\times n}$ endowed with the inner product
$\langle M, N \rangle = \trace{M^{T}N} = \sum_{i,j=1}^n M_{ij} N_{ij}$
and the corresponding (Frobenius) norm $\| M \|_F = {(\langle M, M \rangle)}^{1/2} = {(\sum_{i,j=1}^n M_{ij}^{2} )}^{1/2}$.
Furthermore, given real symmetric matrices $M, N$, we write $M \succeq N$ if and only if $M -N$ is 
positive semi-definite, and $M \succ N$ if and only if $M -N$ is positive definite. 
The spectrum of a symmetric matrix $M$ is denoted by $\sigma(M)$.

\section{The LQR Problem and the Riccati Equation} 
\label{sec:riccati}

In this section we present the mathematical statement of the LQR problem and the governing 
Hessenberg index-2 DAE, and we show
how it can be transformed into a `standard' LQR problem governed by an ODE using a projection onto the 
constraint manifold of the original DAE. Then we apply classical LQR theory to this transformed
problem to compute, under standard conditions on the system, the solution of the LQR problem
via the GCARE. As mentioned before, the problem transformation is performed to derive the 
solution, but the computations are done using the original DAE framework. 
The projection used to convert the DAE into an ODE was first used in a different
context by Heinkenschloss et al.\ \cite{morHeiSS08}. For DAEs derived from a finite 
element discretization of the Stokes or linearized Navier-Stokes system, 
B{\"a}nsch et al.\ \cite{BaeBSetal15} show that this projection is a discrete version of the Leray projector. 
Projections have also been used by Benner, Stykel \cite{BenS14} to formulate and solve GCAREs associated with index-2 DAEs, 
although, as already noted in the introduction, the projection there is different.
Except for some extensions in problem statement and notation the
material in this section is mostly known from   \cite{BaeBSetal15,morHeiSS08}, 
but is needed  to provide the necessary background that allows us to switch between expressions
using the original DAE system and the the corresponding expressions using the transformed ODE system.
Compared to \cite{BaeBSetal15}, this section also provides a more detailed link between the representations
of the optimal control of the LQR problem derived using the original DAE and transformed ODE system.

\subsection{The LQR Problem}     \label{sec:LQR}
Given matrices $A, M  \in \R^{n_v\times n_v}$,
$G \in \R^{n_v\times n_p}$,
$B \in \R^{n_v\times n_u}$, and
$C \in \R^{n_y\times n_v}$ such that
$M$ is symmetric positive definite and $G$ has rank $n_p < n_v$,
we consider the LQR problem
\begin{equation}   \label{eq:LQR-objective}
        \min_{\bu \in L^2(0,\infty)} \int_0^\infty  \|  \by(t) \|_2^2 +  \|  \bu(t) \|_2^2 \, dt,
\end{equation}
where for given $\bu \in L^2(0,\infty)$, the function $\by \in L^2(0,\infty)$ is obtained as the 
output of the  Hessenberg index-2 Differential Algebratic Equation system
\begin{subequations}\label{eq:DAE}
\begin{align}
    M\frac{d}{dt}\bv(t)&=A\bv(t)+G\bp(t)+B \bu(t),    \label{eq:DAE1}\\
                              0&=G^T \bv(t),                      \label{eq:DAE2}\\
                        \by(t)&= C\bv(t).                       \label{eq:DAE3}
  \end{align}
\end{subequations}
To ensure well-posedness of the LQR, we will make additional assumptions on the 
system (\ref{eq:DAE}) in Section~\ref{sec:LQR-sol}. In the cost functional, we may replace the Euclidian norms by any
weighted norm induced by positive definite matrices \(Q_y\) and \(Q_u\). Here, we set both weighting matrices
to the appropriate identity for ease of notation. It is straightforward to include non-identify weighting
matrices into the problem description and the computational framework.

The  LQR problem (\ref{eq:LQR-objective}, \ref{eq:DAE}) arises, e.g., in feedback stabilization of the 
Navier-Stokes equations, see B{\"a}nsch et al.\ \cite{BaeBSetal15} or Raymond \cite{Ray06}.
In this context, (\ref{eq:DAE1}, \ref{eq:DAE2}) correspond to the linearized discretized
Navier-Stokes equations, and $\bv$, $\bp$ correspond to  velocity and pressure, respectively.
The problem also arises in feedback stabilization of multi-field flow problems, see
B{\"a}nsch et al.\  \cite{BaeBSetal14}.
In this case, (\ref{eq:DAE1}) includes additional equations such as linearized reaction equations,
and $\bv$ corresponds to velocities and the other fields, such as concentrations. 

If we define 
\begin{equation}  \label{eq:DAE-matrices-compact}
    \bA = \begin{bmatrix} A &G\\G^T &0\end{bmatrix}, \quad
    \bM = \begin{bmatrix} M & 0 \\ 0 &0\end{bmatrix}, \quad
    \bB = \begin{bmatrix} B \\  0 \end{bmatrix}, \quad
    \bC = \begin{bmatrix} C & 0 \end{bmatrix}, 
\end{equation}
and
\[
     \bx(t) = \begin{bmatrix} \bv(t) \\ \bp(t) \end{bmatrix},
\]
the DAE system (\ref{eq:DAE}) can be written in the  compact form
\begin{subequations}\label{eq:DAE-compact}
\begin{align}
    \bM\frac{d}{dt}\bx(t)&=\bA\bx(t) + \bB\bu(t),   \\
                      \by(t)&=\bC\bx(t).       
  \end{align}
\end{subequations}

The structure of  (\ref{eq:DAE}) can be used to convert  the  LQR problem (\ref{eq:LQR-objective}, \ref{eq:DAE}) 
into a classical one governed by an ODE. We proceed as in \cite{morHeiSS08,BaeBSetal15}.
The constraint (\ref{eq:DAE2}) and the variable $\bp$ can be eliminated from 
(\ref{eq:DAE1}, \ref{eq:DAE2}) via the projection
\begin{align}\label{eq:def-Pi}
  \varPi = I_{n_v}-G(G^T M^{-1}G)^{-1}G^T M^{-1} \in \R^{n_v\times n_v}.
\end{align}
The matrix $\varPi$ obeys $\varPi^2  =  \varPi$  and $\varPi M  =  M \varPi^T$, i.e., it is in fact
an $M$-orthogonal projection.
Furthermore, 
\begin{equation}\label{eq:rangePi-nullG}
  \mbox{null}(\varPi^T )= \mbox{range}(M^{-1}G) \quad
  \text{and}\quad
  \mbox{range}(\varPi^T )= \mbox{null}(G^T ),
\end{equation}
which means that 
\[
      0=G^T \bv(t)   \quad\mbox{ if and only if} \quad
      \bv(t) = \varPi^T  \bv(t).
\]
We use the latter property to enforce (\ref{eq:DAE2}) and multiply (\ref{eq:DAE1})
by $\varPi$ to arrive at
\begin{subequations}\label{eq:DAE-projected}
\begin{align}
    \varPi  M \varPi^T   \frac{d}{dt}\bv(t) =& \varPi A\varPi^T   \bv(t) + \varPi B\bu(t),  \\
                 \by(t) =&C \varPi^T   \bv(t).
 \end{align}
\end{subequations}
If needed, the function $\bp$ can be computed from $\bv$, $\bu$ using
\begin{align}  \label{eq:DAE-p}
   \bp(t) =  - (G^T M^{-1}G)^{-1}  G^T M^{-1} A\bv(t) 
                - (G^T M^{-1}G)^{-1}  G^T M^{-1} B   \bu(t) . 
\end{align}
Equation (\ref{eq:DAE-p}) is obtained multiplying (\ref{eq:DAE1})  by $G^T M^{-1}$ and using (\ref{eq:DAE2}).

Since $\varPi M\varPi^T  \in\R^{n_v\times n_v}$ has an $n_p$-dimensional null-space and cannot be inverted,
(\ref{eq:DAE-projected}) is still not an ODE. However, $\varPi^T   \bv(t) \in \R^{n_v}$ is contained
in the $n_v-n_p$ dimensional subspace  $\mbox{range}(\varPi^T )$ 
and we can explicitly express $\varPi^T   \bv(t)$ as an element of this subspace.
This is done using the decomposition
\begin{align}\label{eq:Pi-decomposition}
  \varPi=\Theta_l\Theta_r^T \quad\text{such that}\quad\Theta_l^T \Theta_r=I_{n_v-n_p}
\end{align}
with $\Theta_l,\Theta_r\in\R^{n_v \times (n_v-n_p)}$. 
In particular
\begin{equation}\label{eq:rangePi}
       \mbox{range}(\Theta_r ) = \mbox{range}(\varPi^T ).
\end{equation}
The new variable $\widetilde{\bv}(t)=\Theta_l^T \bv(t) \in \R^{n_v-n_p}$ satisfies
\begin{equation}\label{eq:Theta*vtilde=v}
          \Theta_r \widetilde{\bv}(t)= \Theta_r \Theta_l^T \bv(t) =  \varPi^T \bv(t) =  \bv(t).
\end{equation}

Using the decomposition (\ref{eq:Pi-decomposition}), we define
\begin{equation}\label{eq:Theta-projected-matrices}
  \cM:=\Theta_r^T M\Theta_r,\quad
  \cA:=\Theta_r^T A\Theta_r,\quad
  \cB:=\Theta_r^T B,\quad
  \cC:=C\Theta_r,
\end{equation}
and write the descriptor system \eqref{eq:DAE-projected} as
\begin{subequations}\label{eq:ODE}
  \begin{align}
    \cM \frac{d}{dt} \widetilde{\bv}(t) &= \cA  \widetilde{\bv}(t)+\cB\bu(t),\\
    \by(t) &= \cC  \widetilde{\bv}(t).
  \end{align}
\end{subequations}

The DAE system (\ref{eq:DAE}) is equivalent to system (\ref{eq:ODE}), which is
an ODE system since with \(M\) being symmetric and positive definite so is
\(\cM\), by \(x^{T}\cM x = {(\Theta_{r}x)}^{T}M(\Theta_{r}x)>0\).
Furthermore, the  LQR problem (\ref{eq:LQR-objective},~\ref{eq:DAE}) is equivalent to
the classical  LQR problem (\ref{eq:LQR-objective},~\ref{eq:ODE}). We summarize this result in the
following proposition.

\begin{proposition}\label{prop:equivalence}
   The functions $\bv$, $\bp$ solve (\ref{eq:DAE1},~\ref{eq:DAE2}) if and only if
           $\bv = \Theta_r \widetilde{\bv}$,  $\widetilde{\bv}$ solves (\ref{eq:ODE}a), and (\ref{eq:DAE-p}) holds.
  Moreover, the control $\bu_{*} \in L^2(0,\infty)$ solves
           the  LQR problem (\ref{eq:LQR-objective},~\ref{eq:DAE}) if and only if it solves
           the classical LQR problem (\ref{eq:LQR-objective},~\ref{eq:ODE})
\end{proposition}

The equivalence between  the  LQR problem (\ref{eq:LQR-objective},~\ref{eq:DAE})
and  the classical LQR problem (\ref{eq:LQR-objective},~\ref{eq:ODE}),  
however, is only used theoretically. Even if the matrices $A, \ldots$ in
(\ref{eq:DAE}) are sparse, the projected matrices $\cA, \ldots$ in (\ref{eq:ODE}) are dense.
We will use the equivalence between (\ref{eq:LQR-objective},~\ref{eq:DAE}) and
(\ref{eq:LQR-objective},~\ref{eq:ODE}) to derive our algorithms, but always compute using
the formulation (\ref{eq:LQR-objective},~\ref{eq:DAE}).

%

\subsection{Solution of the LQR Problem and the Riccati Equation}\label{sec:LQR-sol}
If  $(\cA,\cB;\cM)$ is stabilizable (see Definition~\ref{def:stabilizable}) and
$(\cC,\cA;\cM)$ is detectable (see Lemma~\ref{lem:stab/obser-projected-system}),
the classical LQR problem (\ref{eq:LQR-objective},~\ref{eq:ODE}) has a solution
given as the feedback control law 
\begin{align}\label{eq:LQR-projected-u(t)}
      \bu_{*}(t)=-\underbrace{\cB^T \cX^{(*)} \cM}_{\displaystyle =:\cK^T } \widetilde{\bv}(t),
\end{align}
where $\cX^{(*)} = {(\cX^{(*)})}^T \succeq 0\in \R^{(n_v-n_p)\times (n_v-n_p)}$ 
is the unique stabilizing solution of the GCARE
\begin{align}\label{eq:LQR-projected-GCARE}
        \cC^T \cC+\cA^T \cX\cM+\cM\cX\cA-\cM\cX\cB\cB^T \cX\cM = 0.
\end{align}
See, e.g., Lancaster, Rodman~\cite{LanR95}.

The unique stabilizing solution of the GCARE is obtained by applying Newton's
method to find  a root of the quadratic operator
\begin{align}\label{eq:LQR-projected-GCARE-R}
        \cR(\cX) = \cC^T \cC+\cA^T \cX\cM+\cM\cX\cA-\cM\cX\cB\cB^T \cX\cM.
\end{align}
Given an approximate root $\cX^{(k)}$, the new approximation is computed as the solution
of 
\begin{align}  \label{eq:KNM-Step}
      \cR'(\cX^{(k)}) \cX^{(k+1)}=\cR'(\cX^{(k)}) \cX^{(k)} -\cR(\cX^{(k)}).
\end{align}
This method is known as the Kleinman-Newton method.
See the original paper by Kleinman~\cite{Kle68} or the book by Lancaster, Rodman~\cite{LanR95}.

The system (\ref{eq:KNM-Step}) is a Lyapunov equation and for large-scale problems the 
exact Kleinman-Newton method which is defined by (\ref{eq:KNM-Step}) is impractical.
This is particularly true for the Riccati equation (\ref{eq:LQR-projected-GCARE}) which
is obtained from a large-scale DAE by projection. The projected matrices
in (\ref{eq:Theta-projected-matrices}) are not only large-scale, but because of the projections
they are also dense.
To overcome these difficulties, we need to `undo'  the projections in the numerical computations.
We will discuss the details of our solution approach in the next section.
In the remainder of this section we provide basic
relationships between quantities for the projected problem and quantities for
the original problem.

The Kleinman-Newton method applied to the projected Riccati equation (\ref{eq:LQR-projected-GCARE}) 
generates iterates
\[
       0 \preceq \cX^{(k)} \in  \R^{(n_v-n_p)\times (n_v-n_p)}
\]
and corresponding feedback matrices
\begin{equation}    \label{eq:cKk-def}
       (\cK^{(k)} )^T =  \cB^T \cX^{(k)} \cM  \in  \R^{n_u \times (n_v-n_p)}. 
\end{equation}
We want to write the corresponding feedback law 
$-(\cK^{(k)} )^T \widetilde{\bv}(t) = - \cB^T \cX^{(k)} \cM \widetilde{\bv}(t)$ in terms
of the original variable $\bv = \Theta_r \widetilde{\bv}$, see
Proposition~\ref{prop:equivalence}.
If we define
\begin{subequations}\label{eq:feedback-trans}
\begin{equation}
          X^{(k)} =   \Theta_r  \cX^{(k)}  \Theta_r^T  \in  \R^{n_v \times n_v}
\end{equation}
and
\begin{equation}
      {(K^{(k)})}^T =   B^T X^{(k)} M   \in  \R^{n_u \times n_v},
\end{equation}
then (\ref{eq:Theta-projected-matrices}) and $\bv = \Theta_r \widetilde{\bv}$ imply
\begin{equation}
       {(K^{(k)})}^T \Theta_r  =  {(\cK^{(k)})}^T,
\end{equation}
and
\begin{equation}
       -{(\cK^{(k)})}^T \widetilde{\bv}(t) = - \cB^T \cX^{(k)} \cM \widetilde{\bv}(t)
       = - B^T X^{(k)} M \bv(t) =  -{(K^{(k)})}^T \bv(t).
\end{equation}
\end{subequations}

The convergence of the (exact) Kleinman-Newton method can now be expressed
in the unprojected variables and in the context of the (\ref{eq:DAE}). 
First we show that the stability (detectability) of the system (\ref{eq:ODE}) is equivalent to 
the stability (detectability) of the system (\ref{eq:DAE}).

\begin{definition}\label{def:stabilizable}~
  \begin{enumerate}
  \item A matrix pencil $(\bA,\bM)$  is called  stable
    if it is regular and all the finite eigenvalues of $(\bA,\bM)$ lie
    in the open left half-plane.  
  \item Let $\bA, \bB, \bM$ be given by (\ref{eq:Theta-projected-matrices}). The
    triple $(\bA,\bB;\bM)$ is stabilizable if there exists a matrix $K \in\R^{n_v \times n_u}$
    such that all finite eigenvalues of the matrix pencil 
    \begin{align}\label{eq:initial-closed-loop-pencil}
         \left( \begin{bmatrix}A-B K^T &G\\G^T &0\end{bmatrix},
                 \begin{bmatrix}M&0\\0&0\end{bmatrix}\right)
    \end{align}
    are contained in the open left half-plane.
    The triple  $(\bC,\bA;\bM)$ is called  detectable if and only if
    $(\bA^T ,\bC^T ;\bM)$ is stabilizable. 
  \end{enumerate}
\end{definition}

The following result is proven in~\cite[Lemma~4.4]{Wei16}.
\begin{lemma}\label{lem:stab/obser-projected-system}
  The matrix triple $(\cA,\cB;\cM)$ is stabilizable ($(\cC,\cA;\cM)$
  is detectable) if and only if $(\bA,\bB;\bM)$ is stabilizable
  ($(\bC,\bA;\bM)$ is detectable).
\end{lemma}

With these preparations, the following result is an immediate consequence of the
classical Kleinman-Newton convergence result \cite{Kle68}, \cite{LanR95}. 
See  \cite[Thm.~4.5]{Wei16} for a detailed proof.

\begin{theorem}\label{thm:convergence-KNM-for-projected-Systems}
  Assume $(\bA,\bB;\bM)$ is stabilizable and $(\bC,\bA;\bM)$ is  detectable. 
  There exists a maximal symmetric solution
  $\cX^{(*)} \in  \R^{(n_v-n_p)\times (n_v-n_p)}$ of $\cR(\cX)=0$ for which 
  \begin{equation}   \label{eq:stable-pencil*}
         \left(\begin{bmatrix}A-B B^T X^{(*)} M&G\\G^T &0\end{bmatrix},
               \begin{bmatrix}M&0\\0&0\end{bmatrix}\right)
  \end{equation}
  is stable, where  $X^{(*)}=\Theta_r\cX^{(*)}\Theta_r^T$.
  Furthermore, let  $X^{(0)}=\Theta_r\cX^{(0)}\Theta_r^T$ be symmetric and such that
  \[
    \left(\begin{bmatrix}A-BB^T X^{(0)}M&G\\G^T &0\end{bmatrix},
          \begin{bmatrix}M&0\\0&0\end{bmatrix}\right)
  \]
  is stable, then the sequence $\left\{X^{(k)}\right\}_{k=0}^{\infty}$ defined by
  $X^{(k)}:=\Theta_r\cX^{(k)}\Theta_r^T$, \eqref{eq:KNM-Step}  satisfies
  \begin{align*}
      X^{(1)}\succeq X^{(2)}\succeq\dots\succeq X^{(k)} &\succeq 0,   \\
      \lim_{k\to\infty}X^{(k)}   &=X^{(*)},
  \end{align*}
  and there is a constant $\kappa$ such that
  \[
         \|X^{(k+1)}-X^{(*)} \|_{F} \leq \kappa  \|X^{(k)}-X^{(*)} \|_{F}^{2}
         \quad \mbox{ for all } k.
  \]
\end{theorem}

\begin{remark}
  If $X^{(*)}=\Theta_r\cX^{(*)}\Theta_r^T$ is the solution of the Riccati equation specified
  in  Theorem~\ref{thm:convergence-KNM-for-projected-Systems} and
  \[
       (K^{(*)})^T= B^T X^{(*)} M 
  \]
  is the corresponding feedback matrix, then $\bu_{*}(t):=- (K^{(*)})^T \bv(t)$ 
  solves the LQR problem  (\ref{eq:LQR-objective}, \ref{eq:DAE}).
\end{remark}

In principle, the large-scale projected GCARE (\ref{eq:LQR-projected-GCARE}) can be solved using the 
Kleinman-Newton method \cite{Kle68}. However, the size and special structure
of (\ref{eq:LQR-projected-GCARE}) require the inexact solution of the Newton equation, 
a Lyapunov equation, in each step of the Kleinman-Newton method. Moreover, the
explicit us of the large, dense projected matrices $\cM,  \cA,  \cB,  \cC$
 (\ref{eq:Theta-projected-matrices}) must be avoided in computations and the final
 algorithm must operate with the sparse matrices $\bM,  \bA,  \bB,  \bC$
 (\ref{eq:DAE-matrices-compact}) instead.
 To adopt our approach from \cite{BenHSetal16} to efficiently solve the large-scale projected GCARE 
 (\ref{eq:LQR-projected-GCARE}), we first need to review the main components of our approach in \cite{BenHSetal16}.


\section{Inexact Kleinman-Newton for Algebraic Riccati Equations}
\label{sec:newt_kleinm}
 
 Our approach in  \cite{BenHSetal16} is based on an inexact Kleinman-Newton
 method with line search. Although the exact and, under additional conditions,
 inexact Kleinman-Newton method
 converges with step size fixed to one (see, e.g.,  Kleinman \cite{Kle68}
 or Feitzinger et al.\ \cite{FeiHS09}),
 variable step sizes can hugely improve the performance 
 (Benner, Byers \cite{BenB98}, Benner et al.\ \cite{BenHSetal16}).
 We will also observe this in our numerical tests, see Figure~\ref{fig:accel_Nwt}
 in Section~\ref{sec:numerics}. The line search method and analysis in \cite{BenB98}
 are based on exact Lyapunov equation solves, which guarantees that some favorable
 properties of the Kleinman-Newton iterates are automatically preserved.
 Our paper \cite{BenHSetal16} extends line search algorithms and their analyses to inexact solves. 
 An inexact Kleinman-Newton method without line search is analyzed in \cite{FeiHS09}, but
 some assumptions made in \cite{FeiHS09} do not hold when low-rank methods are applied
 to solve the Lyapunov equation iteratively. 
 We extended the inexact Kleinman-Newton method and analysis 
to integrate the efficient low-rank ADI solver in \cite{BenHSetal16}.
This section reviews the main algorithmic components 
of \cite{BenHSetal16} applied to the projected GCARE (\ref{eq:LQR-projected-GCARE}).
The following Section~\ref{sec:newt_kleinm_DAE2}
then carefully exploits the  special structure of the projected GCARE (\ref{eq:LQR-projected-GCARE})
for an efficient numerical realization.

\subsection{Inexact Kleinman-Newton Method}
\label{subsec:derivation-method}
At its core our method is an inexact Newton method applied to the GCARE $\cR(\cX) = 0$,
where $\cR(\cX)$ is the Riccati residual (\ref{eq:LQR-projected-GCARE-R}).
Given an approximate solution $\cXk \in \R^{(n_v-n_p)\times (n_v-n_p)}$ and 
a so-called forcing  parameter $\eta_k \in(0,1)$,
we compute a step $\cSk  \in \R^{(n_v-n_p)\times (n_v-n_p)}$ that satisfies
\begin{align}\label{eq:inex-KNM-convergence-cond.}
  \Nrm{\cR'(\cXk)\cSk+\cR(\cXk)}\leq\eta_{k}\Nrm{\cR(\cXk)}.
\end{align}
Then we compute a step size $\xi_{k} \in (0,1]$ such that the 
sufficient decrease condition 
\begin{align}\label{eq:suffi-decrease-cond}
  \Nrm{\cR \big(\cXk+\xi_{k}\cSk\big)}\leq(1-\xi_{k}\beta)\Nrm{\cR(\cXk)}
\end{align}
is satisfied, where $\beta>0$ is a given parameter.
The new iterate is
\begin{align}    \label{eq:LS-new-iterate}
      \cXkp =   \cXk+\xi_{k}\cSk.
\end{align}

We will discuss below how we compute an $\cSk$ that satisfies
(\ref{eq:inex-KNM-convergence-cond.}).
As we have shown in \cite{BenHSetal16}, 
if the forcing parameters in \eqref{eq:inex-KNM-convergence-cond.} are
limited by 
\[
         \eta_{k}\leq\bar{\eta}<1 \quad \mbox{ and } \quad
         \beta\in(0,1-\bar{\eta}),
\]
then the sufficient decrease condition (\ref{eq:suffi-decrease-cond}) is satisfied 
for all step sizes $\xi_{k}$ 
\begin{align}\label{eq:step-size-limit}
  0<\xi_{k}\leq(1-\bar{\eta}-\beta)\frac{\Nrm{\cR(\cXk)}}{\Nrm{\cM\cSk\cB\cB^{T}\cSk\cM}}.
\end{align}
To ensure convergence of the sequence of iterates $\{ \cXk \}$, 
the step sizes $\xi_{k}$ also need to be bounded away
from zero. We will state the precise convergence result later, see Theorem~\ref{thm:conv} below.
We use the Armijo rule to compute the step sizes $\xi_{k}$. 
This step size rule and others are discussed in \cite{BenHSetal16}, as well as
conditions that ensure $\xi_{k} \ge\xi_{\min}>0$ for all $k$.

Instead of computing the new iterate $\cSk$ as an approximate solution of
$\cR'(\cXk)\cSk = -\cR(\cXk)$, it is more favorable for our purposes to
compute
\begin{align}\label{eq:prelim-sol.}
  \wtcXkp:=\cXk+\cSk
\end{align}
as an approximate solution of
$\cR'(\cXk)  \wtcXkp = -\cR(\cXk) + \cR'(\cXk)  \cXk$. 
Both equations $\cR'(\cXk)\cSk = -\cR(\cXk)$
and $\cR'(\cXk)  \wtcXkp = -\cR(\cXk) + \cR'(\cXk)  \cXk$ are Lyaponov equations, but the
right hand side of the latter equation, 
\[
      -\cR(\cXk) + \cR'(\cXk)  \cXk
      = -\cC^T \cC - \cM\cXk\cB\cB^T \cXk\cM
      = -\begin{bmatrix}\cC^{T}&\cKk\end{bmatrix}
          \begin{bmatrix}\cC^{T}&\cKk\end{bmatrix}^T,
\]
where $\cK^{(k)}$ is defined in (\ref{eq:cKk-def}), is low-rank and this will allow the
application of the efficient low-rank ADI method (discussed in the next section) 
to compute  $\wtcXkp$. Note that
\[
       \cR'(\cXk)  \wtcXkp =  \big(\cAk\big)^T\wtcXkp\cM+\cM\wtcXkp\cAk,
\]
where
\begin{equation}   \label{eq:cAk-def}
     \cAk = \cA - \cB \cB^T \cXk \cM = \cA - \cB \, \big( \cK^{(k)} \big)^T.
\end{equation}

We define the projected Lyapunov residual at any $ \wtcXkp$ by
\begin{align}\label{eq:proj.-Lyap-Res.-Def.}
          \cLwtkp := \cR'(\cXk)  \wtcXkp + \cR(\cXk) - \cR'(\cXk)  \cXk = \cR'(\cXk)\cSk+\cR(\cXk).
\end{align}
The inexactness condition \eqref{eq:inex-KNM-convergence-cond.} means that we
have to compute $\wtcXkp$ with
\begin{align}\label{eq:inex-KNM-Step}
    \big(\cAk\big)^T\wtcXkp\cM+\cM\wtcXkp\cAk
    = -\begin{bmatrix}\cC^{T}&\cKk\end{bmatrix} \begin{bmatrix}\cC^{T}&\cKk\end{bmatrix}^T +\cLwtkp
\end{align}
such that the corresponding projected Lyapunov residual satisfies
\begin{align}\label{eq:inexact-Newt-step-cond.}
  \Nrm{\cLwtkp}\leq\eta_{k}\Nrm{\cR(\cXk)}.
\end{align}
Using the definition \eqref{eq:LQR-projected-GCARE-R}, \eqref{eq:prelim-sol.}, and \eqref{eq:inex-KNM-Step},
the residual of the projected CARE at \eqref{eq:LS-new-iterate} can be written as
\begin{align}\label{eq:CARE-Taylor-exp.-LS}
    \cR(\cXk+\xi_{k}\cSk)&=\cR(\cXk)+\xi_{k}\cR'(\cXk)\cSk+\frac{\xi_{k}^{2}}{2}\cR''(\cXk)(\cSk,\cSk) \nonumber \\
    &=(1-\xi_{k})\cR(\cXk)+\xi_{k}\cLwtkp-\xi_{k}^{2}\cM\cSk\cB\cB^{T}\cSk\cM,
\end{align}
which can be evaluated efficiently for any $\xi_k$, and therefore can be used to
efficiently compute a step size $\xi_k > 0$ that satisfies \eqref{eq:suffi-decrease-cond}.

The  inexact Kleinman-Newton method with line search is summarized in Algorithm \ref{alg:iKN_method-LS}
below.

\begin{algorithm}[!h]
  \caption{Inexact Kleinman-Newton method with line search} 
  \label{alg:iKN_method-LS}
  \begin{algorithmic}[1]
    \REQUIRE $\cA$, $\cM$, $\cB$, $\cC$, $\tolN$, initial stabilizing
    iterate $\cX^{(0)}$, $\bar{\eta}\in(0,1)$, and $\beta\in(0,1-\bar{\eta})$ 
    \ENSURE Approximate unique stabilizing solution $\cX^{(*)}$ of GCARE
    \eqref{eq:LQR-projected-GCARE} 
    \STATE Set $k=0$.  
    \WHILE{$\Nrm{\cR(\cXk)} > \tolN$} 
    \STATE $\cKk=\cM\cXk\cB$
    \STATE Set $\cAk= \cA-\cB\cKkT$.  
    \STATE Select $\eta_{k}\in(0,\bar{\eta}]$.
    \STATE Compute $\wtcXkp$ that solves the inexact Lyapunov equation
                 \begin{subequations}  \label{eq:iKN_method_Lyap_solve}
                  \begin{equation}      \label{eq:iKN_method_Lyap_solve-1}
                          \big( \cAk \big)^T\wtcXkp\cM+\cM\wtcXkp\cAk
                          =-\begin{bmatrix}\cC^{T}&\cKk\end{bmatrix}\begin{bmatrix}\cC^{T}&\cKk\end{bmatrix}^T + \cLwtkp
                  \end{equation}
                   \begin{equation}      \label{eq:iKN_method_Lyap_solve-2}
                                \text{with } \Nrm{\cLwtkp}\le\eta_{k}\Nrm{\cR(\cXk)}.
                   \end{equation}
                   \end{subequations}
                    \label{alg:iKN_method_Lyap_solve}
    \STATE Compute $\xi_{k}\in(0,1)$ such that
                 $\Nrm{\cR\big( (1-\xi_{k})\cXk+\xi_{k}\wtcXkp \big)} \leq (1-\xi_{k}\beta) \Nrm{\cR(\cXk)}$.
    \STATE Set $\cXkp=(1-\xi_{k})\cXk+\xi_{k}\wtcXkp$.
    \STATE $k=k+1$
    \ENDWHILE
    \STATE $\cX^{(*)}=\cX^{(k)}$
  \end{algorithmic}
\end{algorithm}

The following convergence theorem for the iterates generated by Algorithm \ref{alg:iKN_method-LS}
is adopted from \cite[Thm.~10]{BenHSetal16} to match the notation of the projected
Riccati equation (\ref{eq:LQR-projected-GCARE}).

\begin{theorem} \label{thm:conv}
   Let $(\cA, \cB; \cM)$ be stabilizable, let $(\cC, \cA; \cM)$ be detectable and
  assume that for all $k$, there exists a symmetric positive
  semi-definite $\wtcXkp$ such that (\ref{eq:inex-KNM-Step}) and
  (\ref{eq:inexact-Newt-step-cond.}) hold.
  Furthermore, let $\cXk$ be the iterates generated by Algorithm \ref{alg:iKN_method-LS}
  and $\cAk= \cA-\cB (\cM\cXk\cB)^T$.
  
  \begin{enumerate}[label=(\roman*),ref=\thetheorem\,(\roman*)]
  \item\label{thm:conv-I} If the step sizes are bounded away from
    zero, i.e., $\xi_{k} \ge\xi_{\min}>0$ for all $k$, then
    $\Nrm{\cR(\cXk)}\rightarrow 0$.
  \item\label{thm:conv-II} If in addition the
    pencils $(\cAk,\cM)$ are stable for $k\ge k_{0}$, and $\cXk\succeq
    0$ for all $k\ge k_{0}$, then $\cXk\to\cX^{(*)}$, where
    $\cX^{(*)}\succeq 0$ is the unique stabilizing solution of the
    GCARE (\ref{eq:LQR-projected-GCARE}).
  \end{enumerate}
\end{theorem}

\subsection{Improved Low-Rank ADI Method}
\label{subsec:improved-low-rank-ADI}
The main expense in the inexact Kleinman-Newton Algorithm~\ref{alg:iKN_method-LS}
is in Step~\ref{alg:iKN_method_Lyap_solve}.
We apply the real low-rank ADI method, which is detailed in  \cite{BenHSetal16} and in \cite[Sec.~6.3.1]{Wei16}.
This method generates a low-rank approximate solution $ \wtcXkp$  of the Lyapunov equation in factored form.
Compared to the original version of the ADI method \cite{BenLP08,LiW02},
which is also the version used in B{\"a}nsch et al.\ \cite{BaeBSetal15}, 
we use two important modifications of the original ADI method.
The first reorganizes the computation to obtain a
low-rank representation of the Lyapunov residual in the ADI
iterations  \cite{BenKS13a}, and the second exploits the fact that the ADI shifts must
occur either as a real number or as a pair of complex conjugate numbers to write
almost all\footnote{The linear system solve still has a complex coefficient
  matrix and thus the intermediate \(\cVl\) is complex. This can be avoided
  along the lines of \cite[Remark 4.4]{Kue16}, but is not done in our implementation.} matrices
in the ADI iterations as real matrices  \cite{BenKS13a}. 
Most importantly, the improved method generates a real matrices $\cZ$ and
$\widetilde{\cW}_\ell$, each with few columns, such that $\cZ\cZ^{T} = \wtcXkp$ satisfies  
(\ref{eq:iKN_method_Lyap_solve-1}) and the corresponding Lyapunov residual  
$\cLwtkp =  \widetilde{\cW}_\ell \widetilde{\cW}_\ell^T$ obeys (\ref{eq:iKN_method_Lyap_solve}).
We refer to \cite{BenHSetal16} or \cite[Sec.~6.3.1]{Wei16} for details on  the derivation of the real low-rank 
ADI method. The detailed listing of this method is given in Algorithm~\ref{alg:Greal_ADI} below.

\begin{algorithm}[h]
  \caption{Generalized real-valued low-rank residual ADI method}
  \label{alg:Greal_ADI}
  \begin{algorithmic}[1]
    \REQUIRE $\cAk,\cKk,\cC$, shifts $\{q_{i}\}_{i=1}^{\ell}=\conj{\{q_{i}\}_{i=1}^{\ell}}\in\C^{-}$
    \ENSURE $\cZ$ such that $\cZ\cZ^{T} = \wtcXkp$ and
                     $\cLwtkp = \widetilde{\cW}_\ell \widetilde{\cW}_\ell^T$ satisfy (\ref{eq:iKN_method_Lyap_solve})
    \STATE Set $\ell=1$, $\cZ=[\,]$, $\widetilde{\cW}_{0} = \begin{bmatrix}\cC^{T}&\cKk\end{bmatrix}$.
    \WHILE{$\Nrm{\widetilde{\cW}_{\ell-1}^{T}\widetilde{\cW}_{\ell-1}} > \eta_{k}\Nrm{\cR(\cXk)}$} 
    \STATE $\cVl=\left(\cAkT+q_{\ell}\cM\right)^{-1}\widetilde{\cW}_{\ell-1}$   \label{alg:Greal_ADI-linsys}
    \IF{$\Imag{q_{\ell}}=0$}
    \STATE $\widetilde{\cW}_{\ell} = \widetilde{\cW}_{\ell-1} - 2q_{\ell}\cM\cVl$\label{alg:Greal_ADI-W-real}
    \STATE $\wtcVl=\sqrt{-2q_{\ell}}\,\cVl$
    \ELSE
    \STATE
    $\gamma_{\ell}=2\sqrt{-\Real{q_{\ell}}},\quad\delta_{\ell}= \Real{q_{\ell}} / \Imag{q_{\ell}}$
    \STATE $\widetilde{\cW}_{\ell} = \widetilde{\cW}_{\ell-1} +
    \gamma_{\ell}^{2}\cM\left(\Real{\cVl}+\delta_{\ell}\Imag{\cVl}\right)$\label{alg:Greal_ADI-W-imag}
    \STATE
    $\wtcVlp=
    \begin{bmatrix}
      \gamma_{\ell}\left(\Real{\cVl}+\delta_{\ell}\Imag{\cVl}\right)
      &\gamma_{\ell}\sqrt{(\delta_{\ell}^{2}+1)}\Imag{\cVl}\end{bmatrix}$
    \STATE $\ell=\ell+1$
    \ENDIF
    \STATE $\cZ=\begin{bmatrix}\cZ&\wtcVl\end{bmatrix}$
    \STATE $\ell=\ell+1$
    \ENDWHILE
  \end{algorithmic}
\end{algorithm}

Algorithms~\ref{alg:iKN_method-LS} and \ref{alg:Greal_ADI} work with the projected matrices,
but need to be  implemented operating on the matrices $\bM,  \bA,  \bB,  \bC$. This transformation
will be described in the next section.



\section{Inexact Kleinman-Newton for Algebraic Riccati Equations  Associated with Index-2 DAEs}
\label{sec:newt_kleinm_DAE2}

The inexact Kleinman-Newton Algorithm~\ref{alg:iKN_method-LS} and 
the improved ADI Algorithm~\ref{alg:Greal_ADI} are derived and stated
in terms of the  projected matrices in \eqref{eq:Theta-projected-matrices}. 
As stated before, these matrix are dense, expensive to compute with
and the explicit use of the projection needs to be avoided. 
As before, we use calligraphic font, like $\cXk$, to denote projected 
quantities, and roman font, like $\Xk$,  to denote the corresponding quantities 
without projection.

Regarding the transformation of the iterates in the inexact Kleinman-Newton Algorithm~\ref{alg:iKN_method-LS},
we already know from \eqref{eq:feedback-trans} that
\begin{subequations}\label{eq:matrix-trans-KN}
\begin{equation}
          X^{(k)} =   \Theta_r  \cX^{(k)}  \Theta_r^T  \in  \R^{n_v \times n_v},
\end{equation}
\begin{equation}
       (K^{(k)} )^T =  B^T X^{(k)} M   \in  \R^{n_u \times n_v}, 
       \quad \mbox{ and } \quad
       (K^{(k)} )^T  \Theta_r =  (\cK^{(k)} )^T.
\end{equation}
\end{subequations}

To undo the projections, we multiply the Lyapunov equations and the Riccati residuals
from the left by $\Theta_l$ and from the right by $\Theta_l^T$ and replace
Steps 6 and 7 in Algorithm~\ref{alg:iKN_method-LS} by the following.

\noindent
\rule{\textwidth}{0.4pt}
  \begin{algorithmic}[1]
  \makeatletter
  \setcounter{ALG@line}{5}
  \makeatother
     \STATE Compute $\wtcXkp$ that solves the inexact Lyapunov equation
                 \begin{subequations}  \label{eq:iKN_method_Lyap_solve-trans}
                  \begin{align}      \label{eq:iKN_method_Lyap_solve-1-trans}
                         & \Theta_l \big( \cAk \big)^T\wtcXkp\cM \Theta_l^T   + \Theta_l \cM\wtcXkp\cAk \Theta_l^T      \nonumber \\
                         & = -\Theta_l^T \begin{bmatrix}\cC^{T}&\cKk\end{bmatrix}\begin{bmatrix}\cC^{T}&\cKk\end{bmatrix}^T \Theta_l^T
                              + \Theta_l \cLwtkp \Theta_l^T
                  \end{align}
                  with 
                   \begin{equation}      \label{eq:iKN_method_Lyap_solve-2-trans}
                                 \Nrm{ \Theta_l \cLwtkp \Theta_l^T}\le\eta_{k}\Nrm{ \Theta_l \cR(\cXk) \Theta_l^T}.
                   \end{equation}
                   \end{subequations}
    \STATE Compute $\xi_{k}\in(0,1)$ such that
                 $\Nrm{ \Theta_l \cR\big( (1-\xi_{k})\cXk+\xi_{k}\wtcXkp \big) \Theta_l^T} 
                    \leq (1-\xi_{k}\beta) \Nrm{ \Theta_l \cR(\cXk) \Theta_l^T}$.
  \end{algorithmic}
\rule{\textwidth}{0.4pt}
\medskip

For any symmetric matrix $\cS \in \R^{( n_v - n_p) \times ( n_v - n_p)}$, because $\Theta_l \in \R^{n_v \times ( n_v - n_p)}$ has rank $n_v-n_p$,
$\Theta_l \cS \Theta_l^T = 0$ if and only if $\cS = 0$. 
Thus, replacing Steps 6 and 7 in Algorithm~\ref{alg:iKN_method-LS} by the Steps 6 and 7 above
replaces the Frobenius norm $\| \cdot \|_F$ by the weighted Frobenius norm $\| \Theta_l  \cdot \Theta_l^T \|_F$.
While this change in norm influences the iterates (e.g., because the residual norm is changed
when the inexact Lyapunov equation is solved), it does not change the fundamental
convergence behavior.
In particular, Theorem~\ref{thm:conv} remains valid when the weighted Frobenius norm is used.

\medskip
The reason for multiplying by $\Theta_l$ and $\Theta_l^T$ is that the projection $\varPi$
emerges. In fact, 
using \eqref{eq:Theta-projected-matrices}, \eqref{eq:Pi-decomposition}, and \eqref{eq:matrix-trans-KN},
the left hand side in \eqref{eq:iKN_method_Lyap_solve-1-trans} becomes
\begin{equation}      \label{eq:iKN_method_Lyap_solve-1-trans-trans}
             \Theta_l \big( \cAk \big)^T\wtcXkp\cM \Theta_l^T  + \Theta_l \cM\wtcXkp\cAk \Theta_l^T
              =   \varPi \big( \Ak \big)^T\wtXkp M \varPi^T  + \varPi  M\wtXkp\Ak \varPi ^T,
\end{equation}
where
\[
             \Ak= A-B \big( K^{(k)}  \big)^T. 
\]
Although the projection $\varPi$ emerges in  \eqref{eq:iKN_method_Lyap_solve-1-trans-trans},
it will not be computed and used explicitly.
We outline the main ideas in the following subsections.

\subsection{Low-Rank Residual ADI for Index-2 DAE Systems}
\label{subsec:low-rank-residual-index-2}
Recall \eqref{eq:Pi-decomposition} and \eqref{eq:Theta-projected-matrices}.
We have
\begin{equation}   \label{eq:matrix-trans-KN-W}
     \big[ \cC^{T}\; \cK^{(k)} \big]   = \Theta_r^T   \big[ C^{T}\; K^{(k)} \big] .
\end{equation}
To transform the matrices in the improved ADI Algorithm~\ref{alg:Greal_ADI} we set
\begin{equation}     \label{eq:matrix-trans-ADI-W}
        \widetilde{\cW}_{\ell-1}=\Theta_r^T \widetilde{W}_{\ell-1},\ \ell\geq 1 \mbox{ and }
        \widetilde{W}_{0} :=   \big[ C^{T}\; K^{(k)} \big].
\end{equation}
Using \eqref{eq:Theta-projected-matrices} and \eqref{eq:matrix-trans-KN}, the 
linear system in Step~\ref{alg:Greal_ADI-linsys} of Algorithm~\ref{alg:Greal_ADI}
is transformed into
\begin{equation}    \label{eq:matrix-trans-ADI-linsys1}
     \big( \cA^T+q_{\ell}\cM - \cKk \cB^T \big) \cVl
    =  \Theta_r^T \left(A^T+q_{\ell} M - K^{(k)} B^T\right) \Theta_r \cVl
    =  \Theta_r^T \widetilde{W}_{\ell-1} = \widetilde{\cW}_{\ell-1}.
\end{equation}
We define
\begin{subequations}
\begin{equation}     \label{eq:matrix-trans-ADI-V}
       \Vl = \Theta_r \cVl,\quad\ell\geq 1.
\end{equation}
From \eqref{eq:Pi-decomposition} it follows that
\begin{equation}   \label{eq:matrix-trans-ADI-V1}
      \varPi^T \Vl = \Vl,\quad\ell\geq 1.
\end{equation}
\end{subequations}
Finally, multiplying \eqref{eq:matrix-trans-ADI-linsys1} by $\Theta_l$ from the left,
using \eqref{eq:Pi-decomposition}, \eqref{eq:matrix-trans-ADI-V} and 
\eqref{eq:matrix-trans-ADI-V1}, the linear system in Step~\ref{alg:Greal_ADI-linsys} 
of Algorithm~\ref{alg:Greal_ADI} is written as
\begin{equation}    \label{eq:matrix-trans-ADI-linsys2}
      \varPi \left(A^T+q_{\ell} M - K^{(k)} B^T\right)  \varPi^T \Vl 
    =  \varPi \, \widetilde{W}_{\ell-1}.
\end{equation}
As it is shown by Heinkenschloss et al.\ \cite{morHeiSS08} and  B{\"a}nsch et al.\ \cite{BaeBSetal15}
the solution of the projected system \eqref{eq:matrix-trans-ADI-linsys2} is equivalent to the solution
of the $2 \times 2$ block system
\begin{equation}    \label{eq:matrix-trans-ADI-linsys3}
     \begin{bmatrix} A^T+q_{\ell} M - K^{(k)} B^T &G\\G^T &0\end{bmatrix}
       \begin{bmatrix}  \Vl  \\ *\end{bmatrix}
    =   \begin{bmatrix} \widetilde{W}_{\ell-1} \\ 0\end{bmatrix},
\end{equation}
where ``$*$'' indicates that the second block of the solution matrix is not needed.
Finally, since $K^{(k)} B^T$ is dense, the matrix in \eqref{eq:matrix-trans-ADI-linsys3} 
is written as a low-rank perturbation
\[
     \begin{bmatrix} A^T+q_{\ell} M - K^{(k)} B^T &G\\G^T &0\end{bmatrix}
   =  \begin{bmatrix} A^T+q_{\ell} M  &G\\G^T &0\end{bmatrix}
       -  \begin{bmatrix} K^{(k)} \\ 0\end{bmatrix}  \begin{bmatrix} B^T &  0\end{bmatrix}
 \]
and the solution of \eqref{eq:matrix-trans-ADI-linsys3} is computed using the 
Sherman-Morrison-Woodbury formula. See B{\"a}nsch et al.\ \cite{BaeBSetal15}
or Weichelt \cite[p.~67]{Wei16}.

We use  \eqref{eq:Pi-decomposition}, \eqref{eq:matrix-trans-ADI-W} 
to write the projected Lyapunov residual 
\begin{align}\label{eq:proj.-low-rank-Lyap-res}
    \Theta_l \cL(\widetilde{\cX}^{(k+1)}_\ell) \Theta_l^T
     =  \Theta_l   \widetilde{\cW}_\ell \widetilde{\cW}_\ell^{T} \Theta_l^T
     =  \varPi   \widetilde{W}_{\ell}  \widetilde{W}_{\ell}^{T} \varPi^T
     =:\olWl\olWl^{T}.
\end{align}
Rather than computing $ \widetilde{W}_{\ell}$ and then multiplying by $\varPi$,
we can update $\olWl = \varPi   \widetilde{W}_{\ell} $ directly.
In fact, multiplying line~\ref{alg:Greal_ADI-W-real}
in Algorithm~\ref{alg:Greal_ADI} with $\Theta_l$ from the left and using
 \eqref{eq:matrix-trans-ADI-W}, \eqref{eq:matrix-trans-ADI-V} yields
\[
   \varPi  \widetilde{W}_{\ell}
  = \Theta_l \Theta_r^T  \widetilde{W}_{\ell}
  = \Theta_l \Theta_r^T   \widetilde{W}_{\ell-1}-2q_{\ell} \Theta_l \Theta_r^T M \Theta_r\cVl
  = \varPi   \widetilde{W}_{\ell-1} - 2q_{\ell} \varPi M \Vl
  = \varPi   \widetilde{W}_{\ell-1} - 2q_{\ell}  M \Vl,
\]
where in the last step we have used the $M$-orthogonality of $\varPi$, i.e.,
$\varPi M = M \varPi^T$ and \eqref{eq:matrix-trans-ADI-V1}.
Thus, the projected low-rank residual factor can be accumulated via
\begin{align}\label{eq:real-unprojected-low-rank-Lyap-res}
           \olWl&=\olWlm-2q_{\ell}M\Vl
\end{align}
without using any explicit projections. Only the initial right
hand side $W^{(k)}$ needs to be projected to define
\begin{align}
     \olWn:=\varPi  \; \big[ C^{T}\; K^{(k)} \big] .
\end{align}
This one projection at the beginning of the ADI method is computed by
first solving
\[
     \begin{bmatrix}  M  &G\\G^T &0\end{bmatrix}
    \begin{bmatrix} W \\ * \end{bmatrix}
     =  \begin{bmatrix}  \big[ C^{T}\; K^{(k)} \big] \\ 0\end{bmatrix}  
 \]
 (again, ``$*$'' indicates that the second block of the solution is not used)
 and then setting
 \[
         \olWn = M W. 
 \]
See Heinkenschloss et al.\ \cite{morHeiSS08} or Weichelt \cite[Lemma~4.1]{Wei16}.
This projection is less expensive than a single ADI step and does not considerably increase
the overall computation costs. 
Moreover, the right-hand side $\widetilde{W}_{\ell-1}$ in \eqref{eq:matrix-trans-ADI-linsys2},
\eqref{eq:matrix-trans-ADI-linsys3} can be replaced by $\overline{W}_{\ell-1}$,
since 
\[
     \varPi \, \widetilde{W}_{\ell-1}
     =   \varPi \, \varPi \, \widetilde{W}_{\ell-1}
     =  \varPi \,\overline{W}_{\ell-1}.
\]

To incorporate this improved ADI method into
Algorithm~\ref{alg:iKN_method-LS}, some remaining issues, such as the
storage of the Newton step  and the projected Riccati residual,
need to be addressed. This is done in the next subsection, expanding
the statements in \cite[Sec.~5.2]{BenHSetal16}.

\subsection{Low-Rank Riccati Residual for Index-2 DAE systems}
\label{subsec:low-rank-residuals}
The Newton step $S^{(k)} =   \Theta_r  \cS^{(k)}  \Theta_r^T$ is only used in the computation of the step
size $\xi_{k}$, since the inexact Kleinman-Newton step
\eqref{eq:inex-KNM-Step} directly iterates over the preliminary solution
$\widetilde{X}^{(k)} =   \Theta_r  \widetilde{\cX}^{(k)}  \Theta_r^T$. 
Furthermore, $\cSk$ always occurs in products
$\cM\cSk\cB\in\R^{(n_v-n_p)\times n_{r}}$.  Using
(\ref{eq:LS-new-iterate}), \eqref{eq:prelim-sol.}, and the definition
of the feedback matrix in \eqref{eq:LQR-projected-u(t)}, this product
can be written as
\begin{align}\label{eq:delta-K0}
  \cM\cSk\cB=
  \begin{cases}
    \cM\wtcXkp\cB-\cM\cXk\cB=:\wtcKkp-\cKk=:\dwtcKkp,\quad\xi_{k}\neq 1,\\
    \cM\cXkp\cB-\cM\cXk\cB=:\cKkp-\cKk=:\dcKkp,\quad\xi_{k}=1,
  \end{cases}
\end{align}
which characterizes the feedback change corresponding to the
preliminary or definite new iterate $\wtcXkp$ or $\cXkp$. 
Using \eqref{eq:Theta-projected-matrices}, \eqref{eq:matrix-trans-KN}, and
 $S^{(k)} =   \Theta_r  \cS^{(k)}  \Theta_r^T$, \eqref{eq:delta-K0} becomes
 \begin{align}\label{eq:delta-K}
  M S^{(k)} B=
  \begin{cases}
    M \widetilde{X}^{(k+1)}  B- M X^{(k)}  B=: \widetilde{K}^{(k+1)} - K^{(k)}=: \Delta \widetilde{K}^{(k+1)},\quad\xi_{k}\neq 1,\\
    M X^{(k+1)}  B- M X^{(k)}  B=: K^{(k+1)} - K^{(k)} =: \Delta K^{(k+1)},\quad\xi_{k}=1,\\
  \end{cases}
\end{align}
which characterizes the feedback change corresponding to the
preliminary new iterate or $\widetilde{X}^{(k+1)} $ or new iterate  $X^{(k+1)}$. 
Hence, the dense Newton step $S^{(k)}$ is never formed explicitly. 

The definition $\wtcXkp = \cZ\cZ^{T}$ and update in Step~13 of Algorithm~\ref{alg:Greal_ADI}
implies the formula  
\begin{equation}     \label{eq:wXkp1-update}
\wtcXkp_\ell = \wtcXkp_{\ell-1} + \widetilde{\cV}_\ell \widetilde{\cV}_\ell^{T}, \qquad \ell \ge 1,
\end{equation}
for the implicit iterate $\wtcXkp_\ell$ in Algorithm~\ref{alg:Greal_ADI}.
Algorithm~\ref{alg:Greal_ADI} and  \eqref{eq:matrix-trans-ADI-V} lead to the definition
\begin{equation}     \label{eq:matrix-trans-ADI-Vtrans}
       \widetilde{V}_\ell = \Theta_r \widetilde{V}_\ell,\quad\ell\geq 1.
\end{equation}
Finally, \eqref{eq:delta-K}, \eqref{eq:wXkp1-update}, \eqref{eq:matrix-trans-KN}, and \eqref{eq:matrix-trans-ADI-V}
imply that the feedback change can be accumulated during the ADI algorithm as follows
\begin{align}\label{eq:feedback-change-accumulation}
    \dwtKkpl&=\wtKkpl-K^{(k)}=\wtKkplm+M\wtVl(\wtVl^{T}B)-K^{(k)}   \nonumber \\
    &=\dwtKkplm+M\wtVl(\wtVl^{T}B), & \forall\ell\geq 1
\end{align}
with $\dwtKkpn=-K^{(k)}$; compare \cite[Sec.~5.2]{BenHSetal16}. 
If we consider the  feedback change at the final ADI iteration $\ell$, we simply write
$\dwtKkp$ instead of  $\dwtKkpl$.

The  Riccati residual can be written in low-rank form as
\begin{subequations}\label{eq:real-low-rank-Ric-residual-k}
  \begin{align}
    \cR(\cXk)=\cWk\cWkT-\dcKk\dcKkT=:\cUk\cD\cUkT
  \end{align}
 with
  \begin{align}
    \cUkp=\begin{bmatrix}\cWkp&\dcKkp\end{bmatrix},\quad
    \cD=\begin{bmatrix}I&0\\0&-I\end{bmatrix}.
  \end{align}
\end{subequations}
This representation can be used to efficiently compute $\Nrm{ \Theta_l \cR(\cXk) \Theta_l^T}$.

In the initial iteration $k=0$ with $\cX^{(0)} = 0$, \eqref{eq:real-low-rank-Ric-residual-k} 
holds with $\cW^{(0)}= \cC^T$ and $\dcKk =0$.
Equation \eqref{eq:CARE-Taylor-exp.-LS} and 
$\cL(\wtcXkp)=\widetilde{\cW}_\ell\widetilde{\cW}_\ell^T$
imply
\begin{align}\label{eq:low-rank-Ric-res-LS}
    &\cR(\cXkp)=\cR(\cXk+\xi_{k}\cSk)    \nonumber\\
    &=(1-\xi_{k})\cUk\cD\cUkT+\xi_{k}\widetilde{\cW}_\ell\widetilde{\cW}_\ell^T
    -\xi_{k}^{2}\dwtcKkp\dwtcKkpT      \nonumber\\
    &=(1-\xi_{k})\left(\cWk\cWkT-\dcKk\dcKkT\right)+\xi_{k}\widetilde{\cW}_\ell\widetilde{\cW}_\ell^T 
          -\xi_{k}^{2}\dwtcKkp\dwtcKkpT    \nonumber\\
    &=\bigg[
    \Big[\sqrt{(1-\xi_{k})}\,\cWk\;\sqrt{\xi_{k}}\, \widetilde{\cW}_\ell \Big]\;
    \Big[\sqrt{(1-\xi_{k})}\,\dcKk\;\xi_{k}\dwtcKkp\Big]
    \bigg]
    \times\begin{bmatrix}I&0\\0&-I\end{bmatrix}  \nonumber\\
    &\phantom{=}\,\times
    \bigg[
    \Big[\sqrt{(1-\xi_{k})}\,\cWk\;\sqrt{\xi_{k}}\,\widetilde{\cW}_\ell \Big]\;
    \Big[\sqrt{(1-\xi_{k})}\,\dcKk\;\xi_{k}\dwtcKkp\Big]
    \bigg]^{T},
\end{align}
which is of the form \eqref{eq:real-low-rank-Ric-residual-k} with
\begin{align*}  
    \cWkp:=\begin{bmatrix}
                   \sqrt{(1-\xi_{k})}\,\cWk&\sqrt{\xi_{k}}\,\widetilde{\cW}_\ell \vphantom{\dwtcKkp}
                 \end{bmatrix}, \quad
    \dcKkp:=\begin{bmatrix}
                    \sqrt{(1-\xi_{k})}\,\dcKk&\xi_{k}\dwtcKkp
                  \end{bmatrix}.
\end{align*}

Using \eqref{eq:matrix-trans-ADI-W},  \eqref{eq:low-rank-Ric-res-LS} 
the projected Riccati residual $\cR(\Xkp) := \Theta_l \cR(\cXkp) \Theta_l^T$ can be written as
\begin{align}\label{eq:proj.-low-rank-ric-res}
  \cR(\Xkp) & := \Theta_l \cR(\cXkp) \Theta_l^T  \nonumber\\
                 & = \Theta_l \cWkp\cWkpT \Theta_l^T  - \Theta_l \dcKkp\dcKkpT \Theta_l^T   \nonumber\\
                 &=\varPi\Wkp\WkpT\varPi^T-\varPi\dKkp\dKkpT\varPi^T   \nonumber\\
                 &=\olWkp\olWkpT-\dKkp\dKkpT   
                 =:U^{(k+1)}\cD\left(U^{(k+1)}\right)^{T}
\end{align}
with $U^{(k+1)}=\begin{bmatrix}\olWkp&\dKkp\end{bmatrix}$.
In the second to last equation in \eqref{eq:proj.-low-rank-ric-res} we have used the identity
\begin{align}\label{eq:proj-inv.feedb-change}
      \varPi\dKkp=\varPi M(\Xkp-\Xk)B=M\varPi^T(\Xkp-\Xk)B=\dKkp,
\end{align}
which follows from the $M$-orthogonality of $\varPi$ and $\varPi^T(\Xkp-\Xk) = \Xkp-\Xk$
(cf.\ \eqref{eq:matrix-trans-ADI-V1}).
The updates of  $\cWkp$ and $\dcKkp$ imply
\begin{align}
  \begin{aligned}
    \olWkp&:= \begin{bmatrix}  \sqrt{1-\xi_{k}}\;\olWk&\sqrt{\xi_{k}}\;\wtWl  \end{bmatrix},\quad&\xi_{k}\in(0,1],\\[1ex]
    \olW^{(0)}&:= \varPi \begin{bmatrix}C^{T}&K^{(0)}\end{bmatrix}, \\
    \dKkp&:= \begin{bmatrix}   \sqrt{1-\xi_{k}}\dKk&\xi_{k}\dwtKkp  \end{bmatrix},\quad&\xi_{k}\in(0,1],
  \end{aligned}
\end{align}
where $K^{(0)}$ is an initial stabilizing feedback.

\begin{algorithm}[p]
  \setlength{\abovedisplayskip}{.25em}
  \setlength{\belowdisplayskip}{.25em}
  \caption{Inexact low-rank Kleinman-Newton-ADI for index-2
    DAE systems}
  \label{alg:i_lr_KN_ADI}
  \begin{algorithmic}[1]
    \REQUIRE $M,A,G,B,C$, initial feedback $K^{(0)}$, $\tolN$,
    $\bar{\eta}\in(0,1)$, and $\beta\in (0, 1-\bar{\eta})$
    \ENSURE feedback matrix $K$
    \STATE Set $\olW^{(0)}=\varPi \begin{bmatrix} C^{T}&K^{(0)}\end{bmatrix}$, $\Delta K^{(0)}=0$, 
                 $U^{(0)}=  \begin{bmatrix}\olW^{(0)}&\Delta K^{(0)}\end{bmatrix}$.
    \STATE Set $k=0$.
    \WHILE{$\Big(\Nrm{U^{(k)} \cD \left(U^{(k)}\right)^{T}}\vphantom{^{T}}>  \tolN\Nrm{U^{(0)} \cD \left(U^{(0)}\right)^{T}}\Big)$} 
    \STATE \parbox[t]{.95\linewidth}{Compute ADI shifts
    $\{q_{i}\}_{i=1}^{n_{\text{ADI}}}=\overline{\{q_{i}\}_{i=1}^{n_{\text{ADI}}}}\subset\C^-$
    ordered such that complex pairs form consecutive entries and choose
    $\eta_{k}\in (0,\bar{\eta}]$.}
    \STATE Set $\wtWn=\varPi \begin{bmatrix} C^{T}&K^{(k)}\end{bmatrix}$,
    $\Delta\widetilde{K}_{0}=-K^{(k)}$. 
    \STATE Set $\ell=1$.
    \WHILE{$\Big(\Nrm{\wtWlm^{T}\wtWlm}>\eta_{k}\Nrm{U^{(k)} \cD \left(U^{(k)}\right)}\Big)$}
    \STATE Get $\Vl$ by solving
    $$\begin{bmatrix}A^{T}-K^{(k)}B^{T}+q_{\ell}\,M&G\\G^{T}&0\end{bmatrix}
    \begin{bmatrix}
      \Vl\\ *
    \end{bmatrix}=
    \begin{bmatrix}\wtWlm\\0\end{bmatrix}.$$
    \IF{$\Imag{q_{\ell}}=0$}
    \STATE $\wtWl=\wtWlm-2q_{\ell}M\Vl$\
    \STATE $\wtVl=\sqrt{-2q_{\ell}}\Vl$
    \STATE $\dwtKlp=\dwtKlm+M\wtVl(\wtVl^{T}B)$
    \ELSE 
    \STATE
    $\gamma_{\ell}=2\sqrt{-\Real{q_{\ell}}},\quad\delta_{\ell}= \Real{q_{\ell}} / \Imag{q_{\ell}}$
    \STATE $\wtWlp = \wtWlm +
    \gamma_{\ell}^{2}M\left(\Real{\Vl}+\delta_{\ell}\Imag{\Vl}\right)$
    \STATE
    $\wtVlp=
    \begin{bmatrix}
      \gamma_{\ell}\left(\Real{\Vl}+\delta_{\ell}\Imag{\Vl}\right)
      &\gamma_{\ell}\sqrt{(\delta_{\ell}^{2}+1)}\Imag{\Vl}\end{bmatrix}$
    \STATE $\ell=\ell+1$
    \STATE $\dwtKlp=\dwtKlmm+M\wtVl(\wtVl^{T}B)$
    \ENDIF
    \STATE $\wtUlp=\begin{bmatrix}\wtWlp&\dwtKlp\end{bmatrix}$
    \STATE $\ell=\ell+1$
    \ENDWHILE
    \IF{$\Nrm{\wtUl \cD \wtUl^{T}}>(1-\beta)\Nrm{U^{(k)} \cD \left(U^{(k)}\right)^{T}}$}
    \STATE Compute $\xi_{k}\in(0,1)$ using, e.g., the Armijo rule.
    \ELSE
    \STATE $\xi_{k}=1$.
    \ENDIF
    \STATE
    $\olWkp=\begin{bmatrix}\sqrt{1-\xi_{k}}\;\olWk&\sqrt{\xi_{k}}\;\wtWl\end{bmatrix}$
    \STATE
    $\dKkp=\begin{bmatrix}\sqrt{1-\xi_{k}}\dKk&\xi_{k}\dwtKl\end{bmatrix}$
    \STATE $U^{(k+1)}=\begin{bmatrix}\olWkp&\dKkp\end{bmatrix}$
    \STATE $K^{(k+1)}=(1-\xi_{k})K^{(k)}+\xi_{k}\dwtKl$
    \STATE $k=k+1$
    \ENDWHILE
    \STATE $K=K^{(k)}$
  \end{algorithmic}
\end{algorithm}

Equation \eqref{eq:proj.-low-rank-ric-res} shows that
the Riccati residual $ \cR(\Xkp)$ can be computed without any additional
explicit projection.

The representation \eqref{eq:CARE-Taylor-exp.-LS} shows that
$\cR(\Xk+\xi_{k} S^{(k)})$ is a  quartic polynomial with scalar coefficients.
Just as in \cite[Sec.~5]{BenHSetal16} this is used for an efficient implementation
of the line search computation.

The final feedback at the end of the $k$+$1$-st Newton step is defined via
\begin{align}\label{eq:final-feedback-after-LS}
  K^{(k+1)}=(1-\xi_{k})K^{(k)}+\xi_{k}\dwtKkp.
\end{align}
Only the feedback matrix is needed, but if desired the  Riccati iterate can be
computed in low-rank form as follows.
Assuming the previous Riccati iterate is defined via $X^{(k)}=\Zk\big(\Zk\big)^T$
and the preliminary solution is defined via $\wtXkp=\wtZkp \big(\wtZkp\big)^T$,
the new Riccati iterate can be written as
\begin{align}\label{eq:final-lr-solution-factor-after-LS}
  \begin{aligned}
    \Xkp&=(1-\xi_{k})\Xk+\xi_{k}\wtXkp\\
    &=(1-\xi_{k})\Zk\ZkT+\xi_{k}\wtZkp\big(\wtZkp\big)^T\\
    &=\begin{bmatrix}\sqrt{1-\xi_{k}}\;\Zk&\sqrt{\xi_{k}}\;\wtZkp\end{bmatrix}
    \begin{bmatrix}\sqrt{1-\xi_{k}}\;\Zk&\sqrt{\xi_{k}}\;\wtZkp\end{bmatrix}^{T},
  \end{aligned}
\end{align}
whose size depends on the number of ADI steps in the $k$-th and
$k$+$1$-st Newton iteration. 

The entire process of the inexact low-rank KN-ADI method is depicted
in Algorithm~\ref{alg:i_lr_KN_ADI}.

\section{Numerical Experiments} \label{sec:numerics}

We illustrate the benefits of Algorithm~\ref{alg:i_lr_KN_ADI} to solve  the GCARE associated with the solution of
LQR problem (\ref{eq:LQR-objective}, \ref{eq:DAE})  governed by the linearized Navier-Stokes
equation. 
Since our problem set-up is identical to that in the paper 
by  B{\"a}nsch et al.\ \cite{BaeBSetal15} and in Weichelt's PhD Thesis  \cite{Wei16},
we only sketch it here and refer to \cite{BaeBSetal15,Wei16} for details.
Additional numerical results can be found in  \cite{Wei16}.

\begin{figure}[htb]   
\begin{center}
    \caption{Domain on which the linearized Navier-Stokes equation is posed
                  and coarsest (level 1) triangulation}
  \label{fig:kms-domain}
  \includegraphics[width=0.7\textwidth]{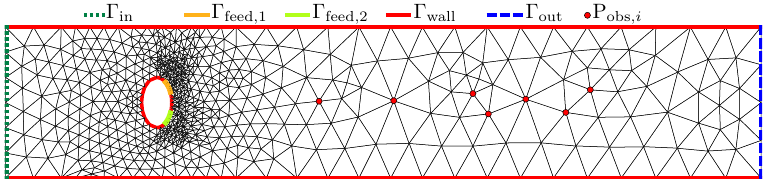}
\end{center}
\end{figure}

The domain on which the Navier-Stokes and linearized Navier-Stokes equations
are posed is shown in Figure~\ref{fig:kms-domain}.
Inflow boundary conditions
are posed on the left boundary, no-slip conditions are posed on part of the cylinder boundary
and on the top and bottom boundary, and outflow conditions are imposed on the right boundary.
Controls are applied 
on two segments on the cylinder wall (indicated by $\Gamma_{{\rm feed},1}$, $\Gamma_{{\rm feed},2}$). Specifically,
for each segment a spatial profile is specified, so that the number of inputs in  (\ref{eq:DAE})
is $n_u = 2$. As described in detail in  \cite[Sec.~2.7]{BaeBSetal15}, \cite[Sec.~4.1.3]{Wei16},
an operator  is constructed that converts these Dirichlet boundary controls to distributed controls, 
such that 
\[
      B \in \R^{n_v \times 2}, \quad n_u = 2
\]
in (\ref{eq:DAE}).
The observations are chosen to be the vertical velocities of the linearized Navier-Stokes equations
at the seven points indicated by $P_{{\rm obs},i}$. Thus, $\by(t) \in \R^7$, $n_y=7$.
Moreover, we penalize the output by $\alpha > 0$, i.e., the output equation (\ref{eq:DAE3})
takes the concrete form
\[
        \by(t) = \alpha C \bv(t)   \quad \mbox{ with }  \quad C \in \R^{7\times n_v}
\]
specified in \cite[p.~A855]{BaeBSetal15}, \cite[Sec.~4.4.1]{Wei16}.

The solution to the steady state Navier-Stokes equation around which is linearized, as well as the 
linearized  Navier-Stokes equations, i.e., the matrices in  (\ref{eq:LQR-objective}, \ref{eq:DAE})
are computed using the finite element flow solver \textsf{NAVIER} \cite{Bae98}, 
which uses  ${\mathcal P}_{2}$--${\mathcal P}_{1}$ Taylor--Hood elements and is written in 
\textsf{FORTRAN90}. The matrices in (\ref{eq:LQR-objective}, \ref{eq:DAE}) are
generated using \textsf{NAVIER} and then stored using the so-called matrix market format \cite{BoiPR96}.
The computations for the resulting matrix equations are performed with
\textsf{MATLAB} R2012b on a 64-bit CentOS~5.5 server with Intel Xeon X5650 at
2.67GHz, with 2~CPUs, 12~cores (6~cores per CPU), and 48~GB main memory available.

We conduct experiments with Reynolds number  $\text{Re} = 100,200,300,400,500$, and 
we use six finite element discretization levels, with Level~1 being the coarsest
(shown in Figure~\ref{fig:kms-domain}). The matrix sizes corresponding to these discretizations
are listed in Table~\ref{tab:refinement_levels}.

\begin{table}[bt]
  \centering
  \caption{Finite element discretization levels and corresponding matrix sizes}
  \label{tab:refinement_levels}
      \begin{tabular}{r|r|r}
      Level & $n_v$ & $n_p$\\ \hline
      1 &  4,796&    672\\
      2 & 12,292& 1,650\\
      3 & 28,914& 3,784\\
      4 & 64,634& 8,318\\
      5 &140,110&17,878\\
      6 &296,888&37,601
    \end{tabular}
\end{table}

For larger Reynolds number, the matrix pencil $(\bA,\bM)$ is not stable
(see \cite[Fig.~2]{BaeBSetal15}, \cite[Sec.~4.2.3]{Wei16})
and a nonzero initial feedback is needed.
We construct the initial feedback $K^{(0)}$ as specified in \cite[Sec.~2.7]{BaeBSetal15}, \cite[Sec.~4.2.3]{Wei16}.

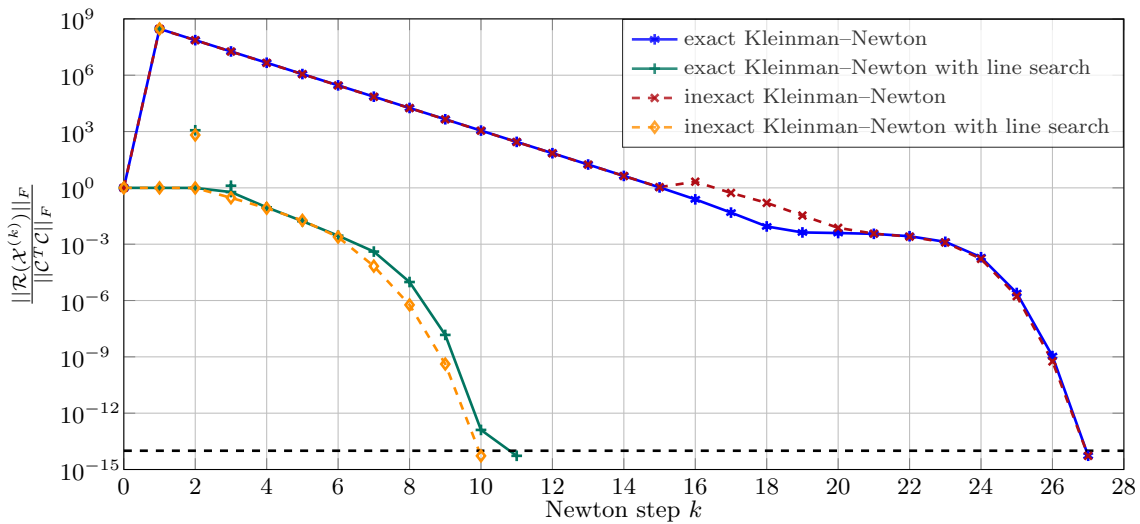
\begin{figure}[tb]   
\begin{center}
    \caption{Impact of the line search on the convergence of the `exact' and 
                  inexact Kleinman-Newton method for the problem with
                  $\text{Re}=500$, Level~1, $\alpha=10^{4},~\tolN=10^{-14}$.
                  There is little difference in the Riccati residuals between the 
                  exact and the inexact Kleinman-Newton method.
                  Line search is active $\xi_k < 1$ for iterations $k=0, 1$ and
                  leads to a dramatic decrease in exact and inexact  Kleinman-Newton iterations.}
  \label{fig:accel_Nwt}
%
%
\begin{tikzpicture}

\definecolor{mycolor1}{rgb}{1,0,1}
\definecolor{tucgreen}{rgb}{0.00392,.458824,.376471}
\definecolor{darkred}{cmyk}{0,1,1,.3}
\definecolor{cscorange}{rgb}{1.0,.5625,0}

\begin{semilogyaxis}[%
scale only axis,
width=.8\linewidth,
height=17em,
ymajorgrids,
xmajorgrids,
xmin=0, xmax=28,
ymin=1e-15, ymax=1e+09,
xlabel={Newton step $k$},
xlabel style={font={\small},at={(.5,0.04)}},
y tick label style={font={\small},anchor=west,xshift=-2.75em},
x tick label style={font={\small}},
legend style={font={\footnotesize},at={(1,1)},nodes=right,anchor=north east,opacity=.8},
ytick={1e-15,1e-12,1e-9,1e-6,1e-3,1e0,1e3,1e6,1e9},
ylabel={$\frac{||\cR(\cX^{(k)})||_{F}}{||\cC^{T}\cC||_{F}}$},
]
\addlegendentry{exact Kleinman--Newton}
\addplot [
color=blue,
solid,line width=1pt,
mark=asterisk,
mark options={solid}
]
coordinates{
 (0,1)(1,2.94655e+08)(2,7.36636e+07)(3,1.84159e+07)(4,4.60398e+06)(5,1.15099e+06)(6,287748)(7,71937)(8,17984.2)(9,4495.97)(10,1123.92)(11,280.908)(12,70.1577)(13,17.4809)(14,4.33281)(15,1.06051)(16,0.2445)(17,0.0482091)(18,0.00881218)(19,0.0042279)(20,0.0039308)(21,0.00356333)(22,0.00265945)(23,0.00133817)(24,0.000201739)(25,2.51641e-06)(26,1.09883e-09)(27,5.41349e-15) 
};

\addlegendentry{exact Kleinman--Newton with line search}
\addplot [
color=tucgreen,
line width=1pt,
mark=+,
mark options={solid}
]
coordinates{
 (0,1)(1,0.999999)(2,0.992382)(3,0.598128)(4,0.0919702)(5,0.0168177)(6,0.00284051)(7,0.000404417)(8,9.67088e-06)(9,1.47385e-08)(10,1.28319e-13)(11,5.36347e-15) 
};
\addlegendentry{inexact Kleinman--Newton}
\addplot [
color=darkred,
solid,line width=1pt,dashed,
mark=x,
mark options={solid}
]
coordinates{
 (0,1)(1,2.91096e+08)(2,7.27741e+07)(3,1.81935e+07)(4,4.54838e+06)(5,1.13709e+06)(6,284274)(7,71068.4)(8,17767.1)(9,4441.73)(10,1110.39)(11,277.545)(12,69.3353)(13,17.3111)(14,4.293)(15,1.09029)(16,2.13307)(17,0.53763)(18,0.158878)(19,0.0335469)(20,0.00742838)(21,0.0036467)(22,0.00252034)(23,0.00123454)(24,0.000165407)(25,1.69798e-06)(26,5.69614e-10)(27,5.33774e-15) 
};
\addlegendentry{inexact Kleinman--Newton with line search}
\addplot [
color=cscorange,
dashed,line width=1pt,
mark=diamond,
mark options={solid}
]
coordinates{
 (0,1)(1,0.999999)(2,0.984461)(3,0.296237)(4,0.0822403)(5,0.0185786)(6,0.00239095)(7,6.77276e-05)(8,5.85334e-07)(9,4.13802e-10)(10,5.33202e-15) 
};
\addplot [
color=tucgreen,
only marks,line width=1pt,
mark=+,
mark options={solid}
]
coordinates{
 (1,2.94655e+08)(2,1168.56)(3,1.30264) 
};
\addplot [
color=cscorange,
only marks,line width=1pt,
mark=diamond,
mark options={solid}
]
coordinates{
 (1,2.91096e+08)(2,653.452) 
};
\addplot [
color=black,line width=1pt,
dashed
]
coordinates{
 (0,1e-14)(28,1e-14) 
};

\end{semilogyaxis}
\end{tikzpicture}

\end{center}
\end{figure}

First, we illustrate the impact of the line search. Figure~\ref{fig:accel_Nwt} shows
the convergence of the `exact' Kleinman-Newton method (i.e., the Lyapunov equation is
solved with fixed high residual tolerance) and the inexact Kleinman-Newton method
(Algorithm~\ref{alg:iKN_method-LS} with $\eta_{k}=\min\{0.1,0.9\cdot \|\cR(X^{(k)})\|_F\}$)
both with and without line search for the LQR problem governed by the discretized
linearized Navier-Stokes equations with  $\text{Re}=500$, output weight $\alpha=10^{4}$,
and discretization level 1. There is little difference in the Riccati  residuals between the 
exact and the inexact Kleinman-Newton method. However, there is a big difference
between the method with and without line search. 
Without line search the relative residual grows dramatically in the initial ($k=0$) iteration. 
With line search, the line search is active $\xi_k < 1$ for iterations $k =0, 1, 2$ (exact Kleinman-Newton)
and iterations $k =0, 1$ (inexact Kleinman-Newton).
Figure~\ref{fig:accel_Nwt} also shows the Riccati residuals corresponding to
$\xi_k = 1$ for the iterations where the line search is active.
That the line search is typically only active in the first few iterations has also been observed in other
applications of Riccati equations (see, e.g., \cite{BenB98}).

Next, we illustrate the influence of the various improvements to the overall 
performance of the Algorithm~\ref{alg:i_lr_KN_ADI}. Specifically we
compare five set-ups, where `Setup~i' corresponds to a basic version
of the Kleinman-Newton-ADI method, and `Setup~v' corresponds to
the most efficient version, which is Algorithm~\ref{alg:i_lr_KN_ADI}. 
Setup~i was used to compute the feedback controls in the paper 
by B{\"a}nsch et al.\ \cite{BaeBSetal15} without explicitly computing the projected residuals (cf., \cite[p. 147]{Wei16}).

The set-ups are given as follows.
\begin{list}{}{}
\item[i:] Kleinman-Newton-ADI method, using the `classic' low-rank ADI
  formulation, with fixed relative 2-norm Lyapunov residual
  tolerance $\tolA=10^{-10}$. ADI shifts are computed heuristically as described
  in \cite{Pen00a}, requiring two short Arnoldi processes to approximate the
  large and small magnitude eigenvalues, i.e. several multiplications and solves
  with the pencil matrices (cf. \cite[Sec.~2.2.3]{Wei16}). Lyapunov
  and Riccati residual norms are computed explicitly (cf. \cite[Sec.~4.3.2]{Wei16}). This algorithm is detailed in  \cite[Sec.~4.2]{Wei16}.
  
\item[ii:] Kleinman-Newton-ADI method with fixed relative 2-norm Lyapunov residual
  tolerance $\tolA=10^{-8}$, real-valued low-rank ADI, heuristic shifts,
  explicit computation of the projected Riccati residual norm.  This explicit
  residual norm computation is not necessary, but demonstrates the accuracy of the
  low-rank Riccati residual.

\item[iii:] Same setup as in ii, except that the low-rank Riccati residual
  updates are used as in \eqref{eq:proj.-low-rank-ric-res}.  The Kleinman-Newton
  and the ADI iterations should be the same in Setup~ii and Setup~iii.

\item[iv:] Same setup as in iii, except that the heuristic shifts are replaced
  by a modified version of the adaptive shifts in \cite{BenKS14b}. At most 15
  ADI shifts are adaptively computed in each call. During the first call (in each
  Newton step), the projected pencil has $n_{y}+n_{u}$ eigenvalues, since we are
  using the right hand side $\wtWn$ for projection. Those eigenvalues are
  passeed into the {\tt lp\_mnmx} routine from \cite{Pen00} to determine
  $r=\min\{15,n_{y}+n_{u}\}$ shifts.

  After all shifts have been used we update the set. To this end,
  the blocks $V_{\ell}$ are stored during the ADI iteration 
  until all previously determined shifts have been used. The entire block
  $Z_{\text{tmp}}=[V_{1},\dots,V_{r}]$ is then used in the adaptive shift
  computation method. A thin QR-decomposition (using {\tt
    qr($Z_{\text{tmp}}$,0)} in \matlab) is performed to determine the new
  projection basis and again upto 15 ADI shifts are determined via {\tt lp\_mnmx}.
\item[v:] Algorithm~\ref{alg:i_lr_KN_ADI} with  $\eta_{k}=\min\{0.1,0.9\cdot\|\cR(X^{(k)})\|_F\}$,
              adaptive shift selection in the ADI method and Armijo line search method.
              (Since the choice $\eta_{k}=\min\{0.1,0.9\cdot\|\cR(X^{(k)})\|_F \}$ of the forcing parameter 
               leads to quadratic convergence of the inexact Kleinman-Newton method
               \cite{BenHSetal16}, this setup
               will also be referred to as `iKNqLS' ({\bf i}nexact {\bf K}leinman {\bf N}ewton 
               with {\bf q}uadratic forcing factor and {\bf L}ine {\bf S}earch).
\end{list}

For each setup, the detailed
iteration numbers (the number of Newton iterations $\#$Newt, 
the number of ADI iterations $\#$ADI, 
and  the number of Newton iterations where the line search was less than one $\#$LS) and the various timings
are depicted in Table~\ref{tab:improvements_timings}. In  Algorithm~\ref{alg:i_lr_KN_ADI}, $k$ is the Newton iteration counter and $\ell$ is the ADI 
iteration counter within a Newton iteration. Note that complex shifts appear as
consecutive pairs for which we solve only one system (in Step~8 of Algorithm~\ref{alg:i_lr_KN_ADI}). 
Still, the ADI iteration counter is increased by two (in Steps~17 and 21 of of Algorithm~\ref{alg:i_lr_KN_ADI}).


  \begin{table}[tb]
    \centering
    \caption{Performance of the various Kleinman-Newton-ADI methods. Iteration
        numbers and timings in seconds for the different Kleinman-Newton-ADI methods
        specified in Setup~i to iv.\ applied to the problems
        with $\text{Re}=500$, refinement level~1, $\tolN=10^{-8}$, and $\alpha=10^{0}$.}
      \label{tab:improvements_timings}
      \begin{tabular}{|c|r|r|r|c|l|l|l|l|}\cmidrule{2-9}
        \multicolumn{1}{c|}{}& $\#\text{KN}$& $\#\text{ADI}$&
        $\#\text{lin\_solve}$&
        $\#\text{LS}$& time$_{\text{lin\_solve}}$
        & time$_{\text{shift}}$ 
        & time$_{\text{proj-res}}$&  {time$_{\text{total}}$}\\\midrule\midrule
       i& 8& 3067& 3067& \multicolumn{1}{c|}{--}& $1.4\cdot10^{3}$& $3.6\cdot10^{1}$
& $5.4\cdot10^{3}$& $6.8\cdot10^{3}$\\\midrule
ii& 8& 3031& 1721& \multicolumn{1}{c|}{--}& $7.0\cdot10^{2}$& $3.6\cdot10^{1}$
& $1.0\cdot10^{1}$&  $7.5\cdot10^{2}$\\\midrule
iii& 8& 3031& 1721& \multicolumn{1}{c|}{--}& $7.0\cdot10^{2}$& $3.7\cdot10^{1}$
&  \multicolumn{1}{c|}{--}&  $7.4\cdot10^{2}$\\\midrule
iv& 8& 600& 346& \multicolumn{1}{c|}{--}& $1.4\cdot10^{2}$& $2.8\cdot10^{0}$
&  \multicolumn{1}{c|}{--}&  $1.5\cdot10^{2}$\\\midrule
v& 7& 305& 176& 1& $7.3\cdot10^{1}$& $1.9\cdot10^{0}$
&  \multicolumn{1}{c|}{--}&  $7.5\cdot10^{1}$\\\midrule
      \end{tabular}
  \end{table}

Comparing Setup~i and Setup~ii in Table~\ref{tab:improvements_timings}
shows that incorporation of the real-valued ADI formulation in Setup~ii
reduces the number of linear solves ($\#$lin\_solve) and, therefore,
the time to solve these systems (time$_{\text{lin\_solve}}$)
drastically. Furthermore, the costs to compute the projected residuals
are reduced by at least two magnitudes. 
Comparing Setup~ii and Setup~iii shows that  avoiding the 
explicit computation of the projected residuals decreases the costs further, 
since the costs to evaluate the low-rank residuals are another magnitude
smaller. 
The adaptive ADI shifts determination in Setup~iv leads to another dramatic
improvement in overall performance. These adaptive shifts reduce the
number of ADI iterations and linear solves by a factor of five. Additionally, the 
computation of these adaptive ADI shifts is one magnitude less expensive
than the heuristic shift computation.

Finally, adding the line-search in Setup~v improves the method further. 
The number of ADI iterations and linear system solves is reduced by a factor of two.
The reduction in ADI iterations also reduced the time for the shift computation.
The line search is less than one only in the first iteration and the cost
of step size computation is negligible. 
Comparing the total computation times shows that the algorithm specified
in Setup~v is 90-times faster than the algorithm specified in Setup~i.
As we will see next, the solution of the Riccati equation becomes more
difficult as the output weighting $\alpha$ increases. In those cases
the speedup of  Setup~v over  Setup~i is even more important.

The following numerical tests focus on  Algorithm~\ref{alg:i_lr_KN_ADI}, 
with Setup~v.  As mentioned before, Algorithm~\ref{alg:i_lr_KN_ADI} with Setup~v
will be referred to as `iKNqLS'.
Table~\ref{tab:iKNADI-steps-NSE} documents the performance of iKNqLS
applied to our test problem as Reynolds number $\text{Re}$, output weight $\alpha$, 
and discretization level change. 
Table~\ref{tab:iKNADI-steps-rey-alpha-NSE}
shows that the number of Kleinman-Newton iterations increases moderately 
with an increasing $\alpha$ and increasing Reynolds number.
Similarly, the number of total ADI iterations needed to approximately solve 
the Lyapunov equations  increases with an increasing $\alpha$ and increasing 
Reynolds number.
Furthermore, line search is only necessary for higher Reynolds numbers and 
higher output weights. 

Table~\ref{tab:iKNADI-steps-rey-Level-NSE} shows that
for $\text{Re} \le 300$, the number of Newton iterations remains
nearly constant as the refinement level is increased.
For $\text{Re} \le 300$ and refinement level greater than two
the number of iterations where the step size is less than one is
unusually large. We believe that this effect is a result of the instability of the matrix pencil. 
Solving the first Newton step inexactly might yield an intermediate solution that is
slightly (especially in finite precision arithmetic) not stabilizing. Therefore, the following ADI iteration tends to diverge. 
Our algorithm detects this behavior by monitoring the Riccati and
Lyapunov residual continuously. 
Although this behavior is not covered by the convergence proof in Theorem~\ref{thm:conv}, 
where a stabilizing solution for $k\geq k_{0}$ is required,  our
algorithm handles this situation by deleting the last ADI step and
performing a line search. This yields convergence in all examples we considered.
The relative Riccati residual seems to stagnate for a couple of iterations and,
hence, an increasing amount of line search runs is required.

\begin{table}[tb]
    \centering
    \caption{Number of Kleinman-Newton iterations {\small ($\#\text{KN}$)}, ADI iterations {\small
        ($\#\text{ADI}$)}, and iterations in which line search was active {\small
        ($\#\text{LS}$)} during the `iKNqLS' process ($\tolN=10^{-8}$,
      $\eta_{k}=\min\{0.1,0.9\cdot\|\cR(X^{(k)})\|_F \}$, Armijo method).}
    \label{tab:iKNADI-steps-NSE}
    \begin{subtable}{\linewidth}
      \centering
      \caption{Influence of output weighting $\alpha$ during
        the `iKNqLS' process (refinement: Level~1).}
      \label{tab:iKNADI-steps-rey-alpha-NSE}
      \begin{tabular}{|l|r|r|r|r|r|r|r|r|r|r|r|r|r|r|r|}\hline
        \multirow{2}{*}{\backslashbox{$\alpha$}{$\text{Re}$}}&\multicolumn{3}{c|}{100$\vphantom{M^{M^{M}}}$}
        &\multicolumn{3}{c|}{200}&\multicolumn{3}{c|}{300}
        &\multicolumn{3}{c|}{400}&\multicolumn{3}{c|}{500}
        \\\cline{2-16}
        &{\tiny$\#\text{KN}$}&{\tiny$\#\text{ADI}$}&{\tiny$\#\text{LS}$}
        &{\tiny$\#\text{KN}$}&{\tiny$\#\text{ADI}$}&{\tiny$\#\text{LS}$}
        &{\tiny$\#\text{KN}$}&{\tiny$\#\text{ADI}$}&{\tiny$\#\text{LS}$}
        &{\tiny$\#\text{KN}$}&{\tiny$\#\text{ADI}$}&{\tiny$\#\text{LS}$}
        &{\tiny$\#\text{KN}$}&{\tiny$\#\text{ADI}$}&{\tiny$\#\text{LS}$}\\\hline
        $\phantom{M}10^{-2\vphantom{M^{M^{M}}}}$ & 3 & 38&\multicolumn{1}{c|}{--}
        & 4 & 74&\multicolumn{1}{c|}{--}
        & 4 & 73 &\multicolumn{1}{c|}{--}
        & 4 & 87 &\multicolumn{1}{c|}{--}
        & 5 & 79 &\multicolumn{1}{c|}{--}\\
        $\phantom{M}10^{-1}$ & 4 & 53&\multicolumn{1}{c|}{--}
        & 5 & 109&\multicolumn{1}{c|}{--}
        & 5 & 84 &\multicolumn{1}{c|}{--}
        & 4 & 74 &\multicolumn{1}{c|}{--}
        & 5 & 109 &\multicolumn{1}{c|}{--}\\\midrule
        $\phantom{M}10^{0}$ & 5 & 80&\multicolumn{1}{c|}{--}
        & 6 & 118&\multicolumn{1}{c|}{--}
        & 7 & 119 &\multicolumn{1}{c|}{--}
        & 6 & 115 & 1
        & 7 & 176 & 1\\\midrule
        $\phantom{M}10^{1}$ & 7 & 98&\multicolumn{1}{c|}{--}
        & 7 & 134&\multicolumn{1}{c|}{--}
        & 8 & 153 & 1
        & 10 & 212 & 2
        & 9 & 201 & 2\\
        $\phantom{M}10^{2}$ & 7 & 109&\multicolumn{1}{c|}{--}
        & 9 & 199& 1
        & 12 & 296 & 3
        & 12 & 331 & 3
        & 12 & 340 & 4\\\midrule
      \end{tabular}
    \end{subtable}\\[1ex]
    \begin{subtable}{\linewidth}
      \centering
      \caption{Influence of refinement levels during the `iKNqLS'
        process ($\alpha=1$).}
      \label{tab:iKNADI-steps-rey-Level-NSE}
      \begin{tabular}{|c|r|r|r|r|r|r|r|r|r|r|r|r|r|r|r|}\hline
        \multirow{2}{*}{\backslashbox{}{$\text{Re}$}}&\multicolumn{3}{c|}{100$\vphantom{M^{M^{M}}}$}
        &\multicolumn{3}{c|}{200}&\multicolumn{3}{c|}{300}
        &\multicolumn{3}{c|}{400}&\multicolumn{3}{c|}{500}
        \\\cline{2-16}
        &{\tiny$\#\text{KN}$}&{\tiny$\#\text{ADI}$}&{\tiny$\#\text{LS}$}
        &{\tiny$\#\text{KN}$}&{\tiny$\#\text{ADI}$}&{\tiny$\#\text{LS}$}
        &{\tiny$\#\text{KN}$}&{\tiny$\#\text{ADI}$}&{\tiny$\#\text{LS}$}
        &{\tiny$\#\text{KN}$}&{\tiny$\#\text{ADI}$}&{\tiny$\#\text{LS}$}
        &{\tiny$\#\text{KN}$}&{\tiny$\#\text{ADI}$}&{\tiny$\#\text{LS}$}\\\hline
        Level 1$\vphantom{M^{M^{M}}}$ & 5 & 80 &\multicolumn{1}{c|}{--}
        & 6 & 118 &\multicolumn{1}{c|}{--}
        & 7 & 119 &\multicolumn{1}{c|}{--}
        & 6 & 115 & 1
        & 7 & 176 & 1\\
        Level 2 & 4 & 73 &\multicolumn{1}{c|}{--}
        & 6 & 118 & 1
        & 7 & 144 & 1
        & 7 & 148 & 1
        & 7 & 168 & 1\\
        Level 3 & 5 & 99 &\multicolumn{1}{c|}{--}
        & 5 & 124 &\multicolumn{1}{c|}{--}
        & 10 & 221 & 3
        & 8 & 200 & 2
        & 7 & 183 &\multicolumn{1}{c|}{--}\\
        Level 4 & 4 & 72 &\multicolumn{1}{c|}{--}
        & 6 & 176 & 1
        & 11 & 198 & 6
        & 10 & 199 & 5
        & 10 & 243 & 3\\
        Level 5 & 5 & 126 &\multicolumn{1}{c|}{--}
        & 6 & 160 & 1
        & 11 & 244 & 4
        & 11 & 273 & 4
        & 10 & 267 & 3\\
        Level 6 & 6 & 189 &\multicolumn{1}{c|}{--}
        & 6 & 184 & 1
        & 11 & 280 & 4
        & 11 & 279 & 4
        & 13 & 344 & 6\\\midrule
      \end{tabular}
    \end{subtable}
  \end{table}

Convergence theory for the exact Kleinman-Newton method
guarantees that the matrix pencils are stable if the initial 
matrix pencil is stable. Thus, another approach to circumvent 
the appearance of a possibly unstable pencil arising from inexact
Lyapunov equation solution is to use a smaller fixed ADI tolerance 
for the first Newton step. 
Rather than using the Lyapunov residual tolerance $\tolA= \eta_k \| \cR(\Xkp) \|$,
we set $\tolA=10^{-2}$ for the first two Newton iterations. 
If the relative Riccati residual decreases and drops below $5\cdot10^{-1}$, the
method switches to the iKNqLS scheme (i.e., $\tolA= \eta_k \| \cR(\Xkp) \|$). 
This is referred to as `exact' start.
Using the $\tolA= \eta_k \| \cR(\Xkp) \|$ in all iterations is referred to as  inexact start.
Table~\ref{tab:iKNADI-start_option} compares  both starting procedures 
for $\text{Re}\geq 300$ and  refinements Level~3--6. 
The `exact' start prevents the stagnation of the relative Riccati
residual and reduces the number of Newton iterations. The line
search is used in at most one iteration.
However, the `exact' solves in the first Newton iterations increase the number of ADI iterations.
Therefore, in most cases a decrease in the  number of Newton iterations does not translate
into a significant decrease in the total number of ADI iterations (and therefore significant
decrease in overall computing time)
when the `exact' start is used.

\begin{table}[tb]
    \centering
    \caption{Comparison of `exact' and
                  inexact start\\ ($\tolN=10^{-8}$,
                   $\eta_{k}=\min\{0.1,0.9\cdot\|\cR(X^{(k)})\|_F \}$, Armijo method, $\alpha = 1$).}
    \label{tab:iKNADI-start_option}
      \begin{tabular}{|c|c|r|r|l|l||r|r|l|l|}\cmidrule{3-10}
        \multicolumn{2}{c|}{}& \multicolumn{4}{c||}{start inexact} &
        \multicolumn{4}{c|}{start "exact" with $\tolA=10^{-2}$}\\\cmidrule{3-10}
        \multicolumn{2}{c|}{}&{$\#\text{KN}$}&{$\#\text{ADI}$}&{$\#\text{LS}$}
        &{time$_{\text{total}}$}
        &{$\#\text{KN}$}&{$\#\text{ADI}$}&{$\#\text{LS}$}
        &{time$_{\text{total}}$} \\\midrule
        \multicolumn{10}{|c|}{$\text{Re}=300$}\\\midrule
        \multirow{6}{*}{\rotatebox{90}{Level}}
        &3& 10& 221& 3&$7.2\cdot10^{2}$
        & 8& 186& 1&$5.9\cdot10^{2}$\\\cmidrule{2-10}
        &4& 11& 198& 6&$1.6\cdot10^{3}$
        & 8& 177& 0&$1.4\cdot10^{3}$\\\cmidrule{2-10}
        &5& 11& 244& 4&$4.8\cdot10^{3}$
        & 8& 215& 0&$4.1\cdot10^{3}$\\\cmidrule{2-10}
        &6& 11& 280& 4&$1.2\cdot10^{4}$
        & 9& 259& 0&$1.2\cdot10^{4}$\\\midrule\midrule
        \multicolumn{10}{|c|}{$\text{Re}=400$}\\\midrule
        \multirow{6}{*}{\rotatebox{90}{Level}}
        &3& 8& 200& 2&$6.1\cdot10^{2}$
        & 6& 158& 1&$5.2\cdot10^{2}$\\\cmidrule{2-10}
        &4& 10& 199& 5&$1.5\cdot10^{3}$
        & 7& 197& 1&$1.6\cdot10^{3}$\\\cmidrule{2-10}
        &5& 11& 273& 4&$5.4\cdot10^{3}$
        & 8& 244& 1&$4.6\cdot10^{3}$\\\cmidrule{2-10}
        &6& 11& 279& 4&$1.3\cdot10^{4}$
        & 8& 272& 1&$1.3\cdot10^{4}$\\\midrule\midrule
        \multicolumn{10}{|c|}{$\text{Re}=500$}\\\midrule
        \multirow{6}{*}{\rotatebox{90}{Level}}
        &3& 7& 183& 0&$6.2\cdot10^{2}$
        & 7& 179& 1&$6.0\cdot10^{2}$\\\cmidrule{2-10}
        &4& 10& 243& 3&$2.0\cdot10^{3}$
        & 8& 192& 1&$1.6\cdot10^{3}$\\\cmidrule{2-10}
        &5& 10& 267& 3&$5.5\cdot10^{3}$
        & 9& 261& 1&$5.5\cdot10^{3}$\\\cmidrule{2-10}
        &6& 13& 344& 6&$1.6\cdot10^{4}$
        & 7& 248& 1&$1.2\cdot10^{4}$\\\midrule\midrule
      \end{tabular}
 \end{table}

Overall iKNqLS is able to solve the Riccati equation in all cases.
Although there is no theoretical justification, our numerics indicate
that the inclusion of line search and computationally inexpensive
monitoring of the low-rank Riccati and Lyapunov residuals enables
the algorithm to successfully cope with intermediate iterates that are
nearly not stabilizing.


\section{Conclusions}
We have extended our inexact Kleinman-Newton method low-rank ADI solver and line search from 
\cite{BenHSetal16}  to  Riccati equations governed by Hessenberg index-2 DAEs. 
Using the projection idea from Heinkenschloss et al.\ \cite{morHeiSS08} and B{\"a}nsch, Benner \cite{BaeB12}
we  transform the  problem governed by the DAE into a `classical'  problem governed  by an ODE.
Our algorithm in \cite{BenHSetal16} is then applied to this transformed problem. However, the
projected ODE is never computed in practice. Instead, a careful exploitation of the problem structure
allows the formulation of the algorithm in the original DAE context. 
We have demonstrated the performance of our Riccati solver to a problem arising in feedback 
stabilization of Navier-Stokes flow around a cylinder. The numerical results document the impact of
various algorithmic components on the overall performance. The algorithmic improvements in this paper 
lead to approximately 90-times speed-up over a previously used Kleinman-Newton-ADI method.
Moreover, we have explored the performance of the new algorithm for various Reynolds numbers,
mesh refinement levels and output weights. The new algorithm was able to solve all instances.
Moreover, although there is no theoretical justification, our numerics indicate
that the inclusion of line search and computationally inexpensive monitoring of the low-rank 
Riccati and Lyapunov residuals enables the new algorithm to successfully cope with intermediate iterates 
that are (slightly) not stabilizing.


\section*{References}
\bibliographystyle{elsarticle-num} 
\bibliography{../../bibfiles/csc,../../bibfiles/mor,../../bibfiles/software}

\end{document}